\documentclass[reqno,11pt]{amsart}
\usepackage{amscd,amssymb,verbatim}

\numberwithin{equation}{section}
\theoremstyle{plain}

\newcommand\br{\langle\, ,\, \rangle}
\newcommand\alp{\alpha}         
\newcommand\bet{\beta}
         \newcommand\Gam{\Gamma}
\newcommand\del{\delta}         \newcommand\Del{\Delta}
\newcommand\eps{\varepsilon}
\newcommand\zet{\zeta}
\newcommand\tet{\theta}         
\newcommand\iot{\iota}

\newcommand\sig{\sigma}         

\newcommand\ome{\omega}         \newcommand\Ome{\Omega}

\newcommand\calA{{\mathcal{A}}}
\newcommand\calB{{\mathcal{B}}}
\newcommand\calC{{\mathcal{C}}}

\newcommand\calE{{\mathcal{E}}}

\newcommand\calH{{\mathcal{H}}}
\newcommand\calI{{\mathcal{I}}}

\newcommand\calM{{\mathcal{M}}}
\newcommand\calN{{\mathcal{N}}}

\newcommand\calP{{\mathcal{P}}}

\newcommand\calR{{\mathcal{R}}}

\newcommand\calT{{\mathcal{T}}}

\newcommand\calX{{\mathcal{X}}}
\newcommand\calY{{\mathcal{Y}}}
\newcommand\calZ{{\mathcal{Z}}}

\newcommand\RR{\mathbb{R}}

\newcommand\ZZ{\mathbb{Z}}

\newcommand\CC{\mathbb{C}}

\newcommand\nek{,\ldots,}
\newcommand\sdp{\times \hskip -0.3em {\raise 0.3ex
\hbox{$\scriptscriptstyle |$}}} 

\newcommand\Cone{\operatorname{Cone}}

\newcommand\Det{\operatorname{Det}}

\newcommand\GL{\operatorname{GL}}

\newcommand\Hom{\operatorname {Hom}}

\newcommand\im{\operatorname {im}}

\newcommand{\Ob}{\operatorname{Ob\,}}

\newcommand\coker{\operatorname{coker}}

\newcommand\Tr{\operatorname{Tr}}

\newcommand\ob{{\overline{b}}}

\newcommand\osig{{\overline{\sig}}}

\newcommand\otet{{\overline{\theta}}}

\newcommand\tilC{{\widetilde{C}}}

\newcommand\tile{{\widetilde{e}}}
\newcommand\tilf{{\widetilde{f}}}

\newcommand\tilK{{\widetilde{K}}}

\newcommand\tilp{{\widetilde{p}}}

\newcommand\tilphi{{\widetilde{\phi}}}

\renewcommand{\>}{\rangle}
\newcommand{\<}{\langle}

\theoremstyle{plain}
\newtheorem{Thm}[subsection]{Theorem}
\newtheorem{Cor}[subsection]{Corollary}
\newtheorem{Lem}[subsection]{Lemma}
\newtheorem{Prop}[subsection]{Proposition}
\newtheorem{Conjec}[subsection]{Conjecture}

\newtheorem{Def}[subsection]{Definition}

\theoremstyle{remark}

\newtheorem{Rem}[subsection]{Remark}

\def\TeXref#1{%
        \leavevmode\vadjust{\setbox0=\hbox{{\tt
                \quad\quad  {\small \textrm #1}}}%
        \theight=\ht0
        \advance\theight by \lineskip
        \kern -\theight \vbox to
        \theight{\rightline{\rlap{\box0}}%
        \vss}%
        }}%
\newif\ifShowLabels
\ShowLabelstrue
\newdimen\theight
\def\TeXrefEq#1{%
        \leavevmode\vadjust{\setbox0=\hbox{{\tt
                \quad\quad  {\small \textrm #1}}}%
        \theight=\ht1
        \advance\theight by \lineskip
        \kern -\theight \vbox to
        \theight{\rightline{\rlap{\box0}}%
        \vss}%
        }}%

\newcommand{\refs}[1]{Section ~\ref{S:#1}}
\newcommand{\refss}[1]{Subsection ~\ref{SS:#1}}
\newcommand{\refsss}[1]{Subsection ~\ref{SSS:#1}}
\newcommand{\reft}[1]{Theorem ~\ref{T:#1}}
\newcommand{\refl}[1]{Lemma ~\ref{L:#1}}
\newcommand{\refp}[1]{Proposition ~\ref{P:#1}}
\newcommand{\refc}[1]{Corollary ~\ref{C:#1}}

\newcommand{\refd}[1]{Definition ~\ref{D:#1}}
\newcommand{\refr}[1]{Remark ~\ref{R:#1}}
\newcommand{\refe}[1]{\eqref{E:#1}}

\newenvironment{thm}[1]%
        { \begin{Thm} \label{T:#1}  \ifShowLabels \TeXref{T:#1} \fi }%
        { \end{Thm} }

\renewcommand{\th}[1]{\begin{thm}{#1}  }
\renewcommand{\eth}{\end{thm} }

\newenvironment{lemma}[1]%
        { \begin{Lem} \label{L:#1}  \ifShowLabels \TeXref{L:#1} \fi }%
        { \end{Lem} }

\newcommand{\lem}[1]{\begin{lemma}{#1} }
\newcommand{\elem}{\end{lemma}}

\newenvironment{propos}[1]%
        { \begin{Prop} \label{P:#1}  \ifShowLabels \TeXref{P:#1} \fi }%
        { \end{Prop} }

\newcommand{\prop}[1]{\begin{propos}{#1} }
\newcommand{\eprop}{\end{propos}}

\newenvironment{corol}[1]%
        { \begin{Cor} \label{C:#1}  \ifShowLabels \TeXref{C:#1} \fi }%
        { \end{Cor} }
\newcommand{\cor}[1]{\begin{corol}{#1}  }
\newcommand{\ecor}{\end{corol}}

\newenvironment{conjec}[1]%
        { \begin{Conjec} \label{Conj:#1}  \ifShowLabels \TeXref{C:#1} \fi }%
        { \end{Conjec} }
\newcommand{\conj}[1]{\begin{conjec}{#1}  }
\newcommand{\econj}{\end{conjec}}

\newenvironment{defeni}[1]%
        { \begin{Def} \label{D:#1}  \ifShowLabels \TeXref{D:#1} \fi }%
        { \end{Def} }
\newcommand{\defe}[1]{\begin{defeni}{#1}  }
\newcommand{\edefe}{\end{defeni}}

\newenvironment{remark}[1]%
        { \begin{Rem} \label{R:#1}  \ifShowLabels \TeXref{R:#1} \fi }%
        { \end{Rem} }
\newcommand{\rem}[1]{\begin{remark}{#1}}
\newcommand{\erem}{\end{remark}}

\newcommand{\eq}[1]%
        { \ifShowLabels \TeXrefEq{E:#1} \fi
           \begin{equation} \label{E:#1} }
\newcommand{\eeq}{\end{equation}}

\newcommand{\prf}{ \begin{proof} }
\newcommand{\eprf}{ \end{proof} }
\newcommand{\Label}[1]{\label{#1}  \ifShowLabels \TeXref{#1} \fi }

\def\Margin#1{%
        \leavevmode\vadjust{\setbox0=\hbox{{\tt
                \quad\quad  {\small \textrm #1}}}%
        \theight=\ht0
        \advance\theight by \lineskip
        \kern -\theight \vbox to
        \theight{\rightline{\rlap{\box0}}%
        \vss}%
        }}%
\newif\ifShowMargin
\ShowMargintrue

\ShowLabelsfalse

\newcommand{\n}{\nabla}

\newcommand{\p}{\text{\( \partial\)}}

\newcommand{\A}{\calA}\newcommand{\B}{\calB}
\renewcommand{\H}{\calH}\newcommand{\E}{\calE}

\renewcommand{\c}{\calC}

\newcommand{\ec}{\calE(\calC)}

\newcommand{\tca}{\calT(\c)}

\renewcommand{\tilp}{\tilde{\p}}

\newcommand{\id}{\operatorname{id}}
\newcommand{\cl}{\operatorname{cl}}

\newcommand{\Hilb}{\text{\it Hilb}}
\renewcommand{\ob}{\operatorname{Ob}}

\begin{document}
\title[$L^2$-torsion and extended $L^2$-cohomology]%
{$L^2$ torsion without the determinant class condition and extended $L^2$
cohomology}
\author[Maxim Braverman]{Maxim Braverman${}^\dag$}
\thanks{${}^\dag$Partially supported by the NSF grant DMS-0204421}
\address{Department of Mathematics\\
        Northeastern University   \\
        Boston, MA 02115 \\
        USA
         }
\email{maxim@neu.edu}

\author{Alan Carey\,${}^\ast$}
\address{Mathematical Science Institute\\ Australian National University\\
Canberra, ACT 0200\\ Australia }
 \email{acarey@intermute.anu.edu.au}

\author[Michael Farber]{Michael Farber${}^\ddag$}
\address{School of Mathematical Sciences\\
Tel-Aviv University\\ Ramat-Aviv 69978, Israel\\
\newline \& Department of Mathematical Sciences, University of Durham, Durham DH1 3LE, UK}
\email{michael.farber@durham.ac.uk}

\thanks{${}^\ddag$Partially supported by a grant from the Israel Science Foundation}
\author{Varghese Mathai${}^\ast$}
\address{Department of Pure Mathematics\\ University of Adelaide\\ Adelaide 5005\\
 Australia}
\email{vmathai@maths.adelaide.edu.au}
\thanks{${}^\ast$ Partially supported by the Clay Mathematics Institute and the
Australian Research Council}

\begin{abstract}
We associate determinant lines to objects of the extended abelian category
built out of a von Neumann category with a trace. Using this we suggest
constructions of the combinatorial and the analytic $L^2$ torsions which,
unlike the work of the previous authors, requires no additional assumptions; in
particular we do not impose the determinant class condition. The resulting
torsions are elements of the determinant line of the extended $L^2$ cohomology.
Under the determinant class assumption the $L^2$ torsions of this paper
specialize to the invariants studied in our previous work \cite{CarFarMat97}.
Applying a recent theorem of D. Burghelea, L. Friedlander and T. Kappeler
\cite{BFK01} we obtain a Cheeger - M\"uller type theorem stating the equality
between the combinatorial and the analytic $L^2$ torsions.
\end{abstract}
\maketitle

\section{Introduction}\Label{S:introduction}

M. Atiyah \cite{Atiyah76} was the first who used the concept of von Neumann
dimension in algebraic topology. He applied this notion to measure the size of
the space of harmonic forms on the universal cover. Later S. Novikov and M.
Shubin \cite{NovSh86} suggested more sophisticated spectral invariants which
are now known as the Novikov -- Shubin invariants. These invariants measure the
size of the space of the forms on the universal cover which are only \lq\lq
nearly\rq\rq \, harmonic. Gromov and Shubin \cite{GromSh91} proved the homotopy
invariance of the Novikov -- Shubin invariants.

Around 1986 Novikov and Shubin proposed to study $L^2$ torsion as a von Neumann
analogue of the classical Reidemeister torsion. A detailed treatment of the
theory of $L^2$-torsion was first done in 1992 by A. Carey, V. Mathai and J.
Lott, see \cite{CarMat92}, \cite{Mathai92}, \cite{Lott92}. In these papers the
$L^2$ torsion is a positive real number which is defined under the assumption
that, firstly, the von Neumann Betti numbers vanish and, secondly, certain {\it
\lq\lq determinant class condition\rq\rq} \, is satisfied. The latter condition
can be expressed in term of the spectral density function of the Laplacian; it
is satisfied if the Novikov -- Shubin invariants are positive.

Paper \cite{CarFarMat97} showed how one may modify the construction of the
$L^2$ torsion so that the first of the above mentioned assumptions (the
vanishing of the reduced $L^2$ cohomology) becomes superfluous. Namely, it was
shown in \cite{CarFarMat97} that one may naturally associate determinant lines
to finitely generated Hilbertian modules over von Neumann algebras. Nonzero
elements of these determinant lines are volume forms on Hilbertian modules. The
$L^2$ torsion is then defined in \cite{CarFarMat97} as an element of the
determinant line of the reduced $L^2$ cohomology. However, the determinant
class assumption was still present in \cite{CarFarMat97}.

Let us mention here a standing conjecture (raised by V. Mathai and J. Lott)
that the Novikov -- Shubin invariants are always positive for the special case
of the regular representation $\pi\to \ell^2(\pi)$ of the fundamental group
$\pi$ of a compact polyhedron. At present this conjecture remains open.
However, in the more general setting of Hilbertian representations which we
study in this paper, the determinant class condition may be easily violated
(even when the underlying manifold is the circle). Therefore removing the
determinant class assumption seems to be important for the theory of $L^2$
torsion.

The main goal of this paper is to define the $L^2$ torsion without requiring
the determinant class condition. It turns out that the $L^2$ torsion (both
combinatorial and analytic) is well defined as an element of the determinant
line associated to the extended $L^2$ cohomology. We also prove here the
Cheeger -- M\"uller type theorem stating coincidence between the combinatorial
and the analytic torsions.

The notion of the extended $L^2$ cohomology was first introduced by M. Farber
in \cite{Farber95}, \cite{Farber96}. The main idea was to embed the category of
Hilbertian representation into an abelian category  and to view the cohomology
of the chain complex of $L^2$ chains of Atiyah \cite{Atiyah76} as an object of
the extended abelian category. The work \cite{Farber95}, \cite{Farber96} was
inspired by the general categorical study of P. Freyd \cite{Freyd}. The
approach of \cite{Farber95}, \cite{Farber96} compared with the notion of the
reduced cohomology gives a bigger homological object: the extended cohomology
contains the reduced cohomology and carries the information about the von
Neumann Betti numbers; it also determines the Novikov -- Shubin invariants and
some other new invariants, see \cite{Farbergr}. The great advantage of the
extended cohomology is that it allows to use the well-developed formalism of
abelian categories to compute the $L^2$ cohomology. An algebraic analogue of
the extended cohomology theory of \cite{Farber95}, \cite{Farber96} was later
suggested by W. L\"uck \cite{Luck97}.

The present paper gives an additional reason to believe that the extended
cohomology, allowing to unify the theory of $L^2$ torsion, is the right object
to study. We show how to associate a determinant line $\det \calX$ with any
object $\calX=(\alpha: A'\to A)$ of the extended category. If $\alpha$ is
injective, an element of $\det \calX$ can be represented as a ratio
$$\frac{\br}{\br_1}\, \in \, \det \calX$$
where $\br$ and $\br_1$ are admissible scalar products on $A$ and $A'$
respectively. The determinant line $\det \calX$ is the tensor product of the
determinant lines associated to the projective and torsion parts of $\calX$. We
show in this paper that determinant line of some
 torsion objects (they are called {\it $\tau$-trivial}) can be canonically trivialized.
Moreover, the determinant class condition means precisely that the torsion part
of the extended cohomology is $\tau$-trivial. This fully explains the role of
the determinant class assumption.

Our proof of the Cheeger -- M\"uller type theorem stating the equality between
the analytic and the combinatorial torsions is based on a recent theorem of D.
Burghelea, L. Friedlander and T. Kappeler \cite{BFK01} which uses the language
of the relative torsion.

Since this paper is already quite long we decided to publish the discussion of
examples and the applications elsewhere.

\tableofcontents
\section{Hilbertian von Neumann categories and Fuglede-Kadison
    determinants}\Label{S:hilbcat}

This section contains preliminary material which will be used later in this paper.

First we recall the notions of a von Neumann category and of a trace on a von Neumann
category; here we follow papers \cite{Farber98,Farber00NS}. To understand the main
constructions and results of the paper the reader may keep in mind the simplest and the
most important for applications example of the category of Hilbertian modules over a von
Neumann algebra, cf. \refsss{hilbmod}. We recall also (following \cite{CarFarMat97}) how
one uses a trace to define the Fuglede-Kadison determinant.

\subsection{Hilbertian spaces}\Label{SS:hilb}
Recall that {\em  a Hilbertian space} (cf. \cite{Palais65}) is a topological vector space
$H$, which is isomorphic to a Hilbert space in the category of topological vector spaces.
In other words, there exists a scalar product on $H$, such that $H$ with this scalar
product is a Hilbert space with the originally given topology.

Scalar products with the above property are called {\em  admissible}. Given one
admissible scalar product $\< \ ,\ \>$ on $H$, any other admissible scalar product is
given by
\[
    \<x,y\>_1\ =\ \<Ax,y\>,\quad x,y\in H,
\]
where $A: H\to H$ is an invertible positive operator $A^\ast=A,\quad A>0$.

Hilbertian spaces naturally appear as Sobolev spaces of sections of vector bundles, cf.
\cite{Palais65}.

Let us denote by $\Hilb$ the category of Hilbertian spaces and continuous linear maps.

If $H$ is a Hilbertian space, $H^\ast$ will denote its conjugate, i.e. the space of all
anti-linear continuous functionals on $H$ (bounded $\RR$-linear maps $\phi:H\to \CC$,
such that $\phi(\lambda h)=\overline \lambda\phi(h)$ for all $\lambda\in \CC$ and $h\in
H$; here the bar denotes the complex conjugation). We consider the action of $\CC$ on
$H^\ast$ given by $(\lambda\cdot \phi)(h)=\phi(\overline\lambda\cdot h)$ for all $h\in
H$. The canonical isomorphism
\begin{eqnarray}\Label{E:tag2-2}
    H \ \to \ (H^*)^* \ = \ H^{\ast\ast}
\end{eqnarray}
is given by $h\mapsto(\phi\mapsto\overline{\phi(h)})$, where $h\in H$, and $\phi\in
H^\ast$.

The following definition was suggested in \cite{Farber00NS}.
\defe{Hilbcat} A Hilbertian von Neumann category is an additive subcategory $\c$ of
$\Hilb$ with the following properties:
\begin{enumerate}
\item[(1)]
for any $H\in \ob(\c)$ the dual space $H^\ast$ is also an object of $\c$ and there is a
$\c$-isomorphism $\phi: H\to H^\ast$ such that the formula
\begin{equation}\Label{E:scprH}
    \<x,y\> \ = \ \phi(x)(y), \qquad x, y\in H
\end{equation}
defines an admissible scalar product on $H$;
\item[(2)] for any $H\in \ob(\c)$ the isomorphism \refe{tag2-2} also belongs to $\c$;
\item[(3)] the adjoint of any morphism in $\c$ also belongs to $\c$;
\item[(4)] the kernel $\ker f=\{x\in H; f(x)=0\}$
of any morphism $f: H\to H'$ in $\c$ and the natural inclusion $\ker f \to H$ belong to
$\c$;
\item[(5)] for any $H, H'\in \ob(\c)$ the set of morphisms
\[
    \Hom_{\c}(H, H') \ \subset \ \Hom_{\Hilb}(H, H') \ = \ {\mathcal L}(H, H')
\]
is a linear subspace closed with respect to the weak topology.
\end{enumerate}
\edefe

Thus, objects of $\c$ have structure of Hilbertian spaces and possibly some additional
structure; morphisms of $\c$ are (faithfully) represented by bounded linear maps.

Condition (5) is similar to the well-known condition in the definition of von Neumann
algebras, which explains our term.

\defe{Cadmiss}
Let $H$ be an object of a Hilbertian von Neumann category $\c$. An admissible scalar
product on $H$ is called {\em $\c$-admissible} if it is given by the equation
\refe{scprH} for some $\phi\in \Hom_\c(H,H^*)$.
\edefe
If $\<\ ,\ \>$ and $\<\ ,\ \>_1$ are two $\c$-admissible scalar products then there
exists an invertible positive operator $A\in \Hom_\c(H,H)$ such that
\begin{equation}\Label{E:sc1-sc}
    \<v,w\>_1 \ = \ \<Av,w\>, \qquad\text{for all}\quad v,w\in H.
\end{equation}

Note also, that given any object $H\in\ob(\c)$, a choice of a $\c$-admissible scalar
product on $H$ determines an involution on the algebra $\Hom_{\c}(H,H)$ (adjoint
operator) and the space $\Hom_{\c}(H,H)$ considered with this involution is a von Neumann
algebra.

The conditions (1) - (5) of \refd{Hilbcat} imply the following properties:
\begin{enumerate}
\item[(6)] {\em  The closure of the image $\cl(\im f)$ of any morphism
$f: H\to H'$ in $\c$ and also the natural projection $H'\to H'/\cl(\im f)$ belong to
$\c$.}
\item[(7)] {\em  Suppose that $H'\subset H$ is a closed subspace.
If $H'$, $H$ and the inclusion $H'\to H$ belong to $\c$ then the orthogonal complement
${H'}^\perp$ with respect to a $\c$-admissible scalar product on $H$ and the inclusion
${H'}^\perp\to H$ belong to $\c$.}
\end{enumerate}

\subsection{Finite von Neumann categories}\Label{SS:finitevN}
The following definitions were suggested in \cite{Farber98}.
An object $H$ of a
Hilbertian von Neumann category $\c$ is called {\em finite} if every injective morphism
$f\in \Hom_\c(H,H)$ has a dense image. A Hilbertian von Neumann category is called {\em
finite} iff all its objects are finite.

\subsection{Examples of von Neumann categories}\Label{SS:examples}
\subsubsection{Hilbertian representations}\Label{SSS:hilbrepr}
The simplest and the most important example of a Hilbertian von Neumann category is the
following (see \cite{Farber00NS}, \S 2.3). Let $\mathcal A$ be an algebra over $\CC$ with
involution, which on $\CC$ coincides with the complex conjugation. {\em  A Hilbertian
representation of $\mathcal A$} is a Hilbertian topological vector space $H$ supplied
with a left action $\mathcal A\to {\mathcal L}(H,H)$ of $\mathcal A$ by continuous linear
maps. A morphism $f: H\to H'$ between two Hilbertian representations of $\mathcal A$ is
defined as a bounded linear map commuting with the action of the algebra $\mathcal A$. We
denote the obtained category by $\calC_\calA$.

There is a canonical duality\footnote{A notion of duality in a von Neumann category was
studied in detail in \cite{Farber00NS}.}
in the category of all Hilbertian
representations of a given $\ast$-algebra $\mathcal A$. Namely, given a Hilbertian
representation $H$, consider the space $H^\ast$ of all anti-linear continuous functionals
on $H$. Consider the following action of $\mathcal A$ on $H^\ast$: if $\phi\in H^\ast$
and $\lambda\in \mathcal A$ then $(\lambda\cdot \phi)(h)=\phi(\lambda^\ast\cdot h)$ for
all $h\in H$. Here $\lambda^\ast$ denotes the involution of $\lambda\in \mathcal A$. The
canonical isomorphism
\[
    H \ \to \ H^{**}
\]
is given by $h\mapsto(\phi\mapsto\overline{\phi(h)})$, where $h\in H$, and $\phi\in
H^\ast$. Category of all Hilbertian representations of $\mathcal A$ is a von Neumann
category.


\subsubsection{Hilbertian modules}\Label{SSS:hilbmod}
Let now ${\A}$ be a finite  von Neumann algebra with a fixed finite, normal, and faithful
trace $\tau:{\A}\rightarrow \CC$. In this case there is an important subcategory of the
category $\c_\A$ (described in the previous section).

Let  $*\,$ be the involution in ${\A}$. By $\ell^{2}({\A}) = \ell^{2}_\tau({\A})$ we
denote the completion of ${\A}$ with respect to the scalar product $\<a,b\>=
\tau(b^{*}a)$, for $a,b\in{\A}$, determined by the trace $\tau$. A {\em (projective)
Hilbertian $\A$-module} $M$ is a Hilbertian representation of $\A$ such that there exists
a continuous $\A$-linear embedding of $M$ into $\ell^{2}({\A})\otimes H$, for some
Hilbert space $H$. Note that this embedding is not part of the structure.  A Hilbertian
module $M$ is  said to be {\em finitely generated} if it admits an embedding
$M\rightarrow\ell^{2}({\A})\otimes H$ as above with finite dimensional $H$.

We denote  by $\calH_\calA$ the full subcategory of $\c_\A$, whose objects are Hilbertian
$\A$-modules. We denote by $\calH^f_\A$ the full subcategory of finitely generated
modules in $\calH_\A$.

Note that the algebra $\Hom_{\calH_\A}(M,M)$ coincides with the commutant $\calB=
\B_{\calA}(M)$ of $M$.

$\calH_\A$ is a von Neumann category; $\calH^f_\calA$ is a finite von Neumann category.

\subsubsection{Families of Hilbert spaces}\Label{SSS:familyHs}
Let $Z$ be a locally compact Hausdorff space and let $\mu$ be a positive Radon measure on
$Z$. Let $\mathcal A$ denote the algebra $L^\infty_{\CC}(Z,\mu)$ (the space of
essentially bounded $\mu$-measurable complex valued functions on $Z$, in which two
functions equal locally almost everywhere, are identical). The involution on $\mathcal A$
is given by the complex conjugation. We will construct a category $\c(Z,\mu)$ of Hilbert
representations of $\mathcal A$ as follows. The objects of $\c(Z,\mu)$ are in one-to-one
correspondence with the $\mu$-measurable fields of finite-dimensional Hilbert spaces
$\xi\to\calH(\xi)$ over $(Z,\mu)$, cf. \cite{Dixmier81}, part II, chapter 1. For any such
measurable field of Hilbert spaces, the corresponding Hilbert space is the direct
integral
\begin{equation}\Label{E:bigHs}
    H  \ = \ \int^{\oplus}\, \calH(\xi)\, d\mu(\xi)
\end{equation}
defined as in \cite{Dixmier81}, part II, chapter 1. The algebra $\mathcal A$ acts on the
Hilbert space \refe{bigHs} by pointwise multiplication.

Suppose that we have two $\mu$-measurable finite-dimensional fields of Hilbert spaces
$\xi\to\calH(\xi)$ and $\xi\to\calH'(\xi)$ over $Z$. Then we have two corresponding
Hilbert spaces, $H$ and $H'$, given as direct integrals \refe{bigHs}. We define the set
of morphisms $\Hom_{\c(Z,\mu)}(H,H')$ as the set of all bounded linear maps $H\to H'$
given by {\em decomposable linear maps}
\begin{equation}\Label{E:decompT}
    T \ = \  \int^\oplus\, T(\xi)\, d\mu(\xi),
\end{equation}
where $T(\xi)$ is {\em  an essentially bounded measurable field of linear maps}
$T(\xi):\calH(\xi)\to\calH'(\xi)$, cf. \cite{Dixmier81}, part~II, chapter~2.

It is shown in \cite[\S2.6]{Farber98} that the above construction defines a finite von
Neumann category.

One can also consider a full subcategory $\c^f(Z,\mu)$ of $\c(Z,\mu)$, which consists  of
measurable fields of finite dimensional Hilbert spaces $\H(\xi)$ over $Z$ such that the
dimensions of $\H(\xi)$ are essentially bounded.

Other examples of von Neumann categories can be found in \cite[\S2]{Farber98}.

\subsection{Traces on von Neumann categories}\Label{SS:traces}
Let $\c$ be a Hilbertian von Neumann category.

The following definition was suggested in \S 5.8 of \cite{Farber98}:

\defe{trace}
A trace on $\c$ is a function, denoted $\tau$, which assigns to each object
$\H\in\ob(\c)$ a linear functional
$$\tau_{\H}: \Hom_{\c}(\H,\H)\to \CC$$
on the
von Neumann algebra $ \Hom_{\c}(\H,\H)$,
such that for any pair of objects $\H_1, \H_2\in \Ob(\c)$ the corresponding traces
$\tau_{\H_1}$ and $\tau_{\H_2}$ are compatible in the following sense: if $f\in \hom_{\c}(\H_1,\H_2)$ and
$g\in \hom_{\c}(\H_2,\H_1)$ then
$$\tau_{\H_1} (gf) = \tau_{\H_2}(fg).$$
The trace $\tau$ is called non-negative if $\tau_\H(e)$ is real and non-negative for any idempotent $e\in \hom_\c(\H,\H)$, $e^2=e$.
One says that that a trace $\tau$ on a von Neumann category is normal (or faithful) iff for
each non-zero $\H\in\ob(\c)$ the trace $\tau_\H$ on the von Neumann algebra
$\Hom_{\c}(\H,\H)$ is normal (faithful). \edefe

If $f\in\Hom_{\c}(\H_1\oplus \H_2, \H_1\oplus \H_2)$ is given by a $2\times2$ matrix $(f_{ij})$, where $f_{ij}:\H_i\to\H_j$, $i,j=1,2$, then
\begin{equation}\Label{E:Trofmatr}
    \tau_{\H_1\oplus\H_2}(f)\ =\ \tau_{\H_1}(f_{11}) + \tau_{\H_2}(f_{22})
\end{equation}
as follows easily from the above definition.

To simplify the notation we will often write $\tau$ for $\tau_H$.

\subsection{Example: Trace on the category of finitely generated
        Hilbertian modules}\Label{SS:trfgmod}
Let ${\A}$ be a finite  von Neumann algebra with a fixed finite, normal, and faithful
trace $\tau:{\A}\rightarrow \CC$. Consider the category $\calH^f_\A$ of finitely
generated $\A$-Hilbertian modules, cf. Section~\ref{SSS:hilbmod}.  There is a canonical
trace on this category, which we will also denote by $\tau$, cf.
\cite[Proposition~1.8]{CarFarMat97}. Let us briefly recall the construction of this
trace.

Suppose first that $M$ is {\em free}, that is, $M$ is isomorphic to
$\ell^2(\A)\otimes\CC^k$ for some $k$. Then $\Hom_{\H_A^f}(M,M)$ can be identified with
the algebra of $k\times k$-matrices with entries in $\A$, acting from the right on
$\ell^2(\A)\otimes\CC^k$ (the last module is viewed as the set of row-vectors with
entries in $\ell^2(\A)$). If $\alpha\in\B$ is represented by a $k\times k$ matrix
$(\alpha_{ij}),$ then one define
\[
     \tau_M(\alpha)\ =\ \sum_{i=1}^k \ \tau(\alpha_{ii}).
\]
This gives a trace on $\Hom_{\H_A^f}(M,M)$ which satisfies all necessary conditions.

If $M$ is not free, then we can embed it in a free module as a closed $\A$-invariant
subspace. Then $\Hom_{\H_A^f}(M,M)$ can be identified with a left ideal in the
$k\times{}k$-matrix algebra with entries in $\A$ and the trace described in the previous
paragraph restricts to this ideal and determines a trace on $\Hom_{\H_A^f}(M,M)$. One
then shows that the obtained trace does not depend on the embedding of $M$ in a free
module.

\subsection{Example: Trace on the category of families of Hilbert spaces}\Label{SS:trfamily}
Let $Z$ be a locally compact Hausdorff space endowed with a positive Radon measure $\mu$.
Consider the von Neumann category $\c^f(Z,\mu)$  of measurable fields of finite
dimensional Hilbert spaces $\H(\xi)$ over $Z$ such that the dimensions of $\H(\xi)$ are
essentially bounded, cf. Subsection~\ref{SSS:familyHs}. This category has traces. Let
$\nu$ be a positive measure on $Z$ which is absolutely continuous with respect to $\mu$
and such that $\nu(Z)<\infty$. Then $\nu$ determines the following trace on this von
Neumann category
\[
    \tau_{\nu}(T)\  = \ \int_Z\, \Tr(T(\xi))\, d\nu.
\]
where $T$ is a morphism given by formula \refe{decompT} and $\Tr$ denotes the usual
finite dimensional trace.

\subsection{The Fuglede-Kadison determinant}\Label{SS:Fkdet}
Suppose $\c$ is a Hilbertian von Neumann category endowed with a non-negative, normal,
and faithful trace. Let $M$ be an object in $\c$. Denote by $\GL(M)$ the group of all
invertible elements of $\Hom_\c(M,M)$. We will consider the norm topology on $\GL(M)$;
with this topology it is a Banach Lie group. Its Lie algebra can be identified with
$\Hom_\c(M,M)$. The trace $\tau_{M}$ is a homomorphism of the Lie algebra $\Hom_\c(M,M)$
into the abelian Lie algebra $\CC$. By the standard theorems, it defines a group
homomorphism of the universal covering group of $\GL(M)$ into $\CC$. This approach leads
to following construction of the Fuglede-Kadison determinants, cf.
\cite{HarpeSkandalis84}.

Given an invertible operator $A\in\GL(M)$, find a continuous piecewise smooth path
$A_t\in\GL(M)$ with $t\in [0,1]$, such that $A_0=I$ and $A_1=A$ (it is well known that
the group $\GL(M)$ is pathwise connected, cf. \cite{Dixmier81}). Then define
\begin{equation}\Label{E:logDet}
    \log \Det_{\tau}(A)\ =\ \int_0^1\Re\tau\, [\, A_t^{-1}A_t^\prime\, ]\, dt,
\end{equation}
where $\Re$ denotes the real part.

It has been shown in \cite{CarFarMat97} that the integral does not depend on the choice
of the path, joining $A$ with the identity $I$, and that one has the following
\footnote{In \cite{CarFarMat97} the theorem is stated only for the category of
finitely generated Hilbertian representation of a von Neumann algebra, cf.
\refsss{hilbrepr}. However the proof works for arbitrary von Hilbertian Neumann category
with trace without any changes.} theorem:

\th{FKdet}
The function above,
\begin{equation}\Label{E:Fkdet}
    \Det_{\tau}:\,  \GL(M)\ \longrightarrow \ \RR^{>0}
\end{equation}
called the Fuglede--Kadison determinant, satisfies:
\begin{enumerate}
\item[(a)] $\Det_{\tau}$ is a group homomorphism, that is,
\begin{equation}\Label{E:dethom}
     \Det_{\tau}(A B)\ = \  \Det_{\tau}(A)\cdot  \Det_{\tau}(B)
     \qquad \text{for}\quad A,\ B\in\GL(M);
\end{equation}
\item[(b)]
\begin{equation}\Label{E:detlamI}
    \Det_{\tau}(\lambda I)\ =\ |\lambda|^{\dim_{\A}(M)}
    \qquad\text{for}\quad \lambda \in \CC, \lambda\ne 0;
\end{equation}
here $I\in\GL(M)$ denotes the identity operator;
\item[(c)]
\begin{equation}\Label{E:detlamtau}
     \Det_{\lambda\tau}(A) \ =\ \Det_{\tau}(A)^\lambda
     \qquad\text{for}\quad \lambda\in \RR^{>0};
\end{equation}
\item[(d)]  $\Det_{\tau}$ is continuous as a map $\GL(M)\to \RR^{>0}$, where
$\GL(M)$ is supplied with the norm topology;
\item[(e)] If $A_t$ for $t\in [0,1]$ is a continuous piecewise smooth path in
$\GL(M)$ then
\begin{equation}\Label{E:detA1A0}
     \log \left[\, \frac{\Det_{\tau}(A_1)}{\Det_{\tau}(A_0)}\, \right]\ =\
            \int_0^1\, \Re\Tr\lbrack A_t^{-1}A_t^\prime\rbrack\,  dt.
\end{equation}
Here $A_t^\prime$ denotes the derivative of $A_t$ with respect to $t$.
\item[(f)]  Let $M$ and $N$ be two objects of $\c$, and
$A\in\GL(M)$ and $B\in\GL(N)$ two invertible automorphisms, and let $\gamma:N\to M$ be a
homomorphism. Then the map given by the matrix
\[
   \begin{pmatrix}
         A&\gamma\\ 0&B
   \end{pmatrix}
\]
belongs to $\GL(M\oplus N)$ and
\begin{equation}\Label{E:detM+N}
    \Det_{\tau} \begin{pmatrix} A&\gamma\\ 0&B \end{pmatrix}
    \ =\ \Det_{\tau}(A)\cdot\Det_{\tau}(B).
\end{equation}
\end{enumerate}
\eth

As an example consider an object $M\in \ob(\calC)$ and a $\calC$-automorphism
$A:M\to M$. Fix an admissible scalar product on $M$. Assume that $A$ is
positive and self-adjoint with respect to this scalar product. Let
$A=\int_0^\infty \lambda dE_\lambda$ be the spectral decomposition of $A$. Then
\begin{eqnarray}\label{exdet}
\Det_\tau(A) = \exp\left[\int_0^\infty\ln(\lambda)d\phi(\lambda)\right]
\end{eqnarray}
where $\phi(\lambda)=\tau(E_\lambda)$ is the spectral density function of $A$.

\section{Determinant line of a Hilbertian module}\Label{S:detofmod}

In this section we present a construction of the determinant line associated to
an object in a Hilbertian von Neumann category with trace, which essentially
repeats a construction given in \cite{CarFarMat97}. It will be used later in
\refs{detefc} in the constructions of the determinant line assocoated to an
object in the extended category, and then in Sections~\ref{S:combtor} and
\ref{S:DeRham} to define the combinatorial and the analytic $L^2$-torsions
invariants associated to polyhedra and compact manifolds respectively.

\subsection{The definition}\Label{SS:equiv}
Let $\c$ be a  Hilbertian von Neumann category endowed with a non-negative, normal, and
faithful trace, cf. Definitions~\ref{D:Hilbcat} and \ref{D:trace}
\footnote{We emphasize that for most of our applications it is enough to consider the
special case of the category $\c= \calH_\A^f$ of finitely generated Hilbertian modules,
cf. Subsection~\ref{SSS:hilbmod}.}
. One associates (in a canonical way) with any object $M\in \Ob(\c)$ an oriented real
line which is denoted $\det{M}$; it is the {\em determinant line of} $M$. This
construction generalizes the determinant line of a finite dimensional vector space.

Define $\det{M}$ as a real vector space generated by symbols $\<\ ,\ \>$, one for any
$\c$-admissible scalar product on $M$, subject to the following relations: for any pair
$\<\ ,\ \>_1$ and $\<\ ,\ \>_2$ of $\c$-admissible scalar products on $M$ one writes the
following relation
\begin{equation}\Label{E:equiv}
    \<\ ,\ \>_2  \ = \ {\Det_{\tau}(A)}^{-1/2} \cdot \<\ ,\ \>_1
\end{equation}
where $A\in \GL(M)$ is such that
\[
    \<v,w\>_2\ = \ \<Av,w\>_1
\]
for all $v,w\in M$. Here the transition operator $A\in \Hom_\c(M,M)$ is invertible, cf.
\refe{sc1-sc} and $\Det_{\tau}(A)$ denotes the Fuglede-Kadison determinant of $A$
constructed in \refss{Fkdet} with the aid of the trace $\tau$ on the von Neumann category
$\c$.

Then, cf. \cite[Section~2.1]{CarFarMat97}, $\det{M}$ {\em is the
one-dimensional real vector space generated by the symbol $\<\ ,\ \>$ of any
$\c$-admissible scalar product on $M$}.

Note also, that the real line $\det{M}$ has {\em the canonical orientation}, since the
transition coefficients ${\Det_{\tau}(A)}^{-1/2}$ are always positive. Thus we may speak
of {\em positive and negative} elements of $\det{M}$. The set of all positive elements of
$\det{M}$ will be denoted $\det_+(M)$.

We will think of elements of $\det{M}$ as ``densities" on $M$. If $M$ is a
trivial object, $M=0$, then we set $\det{M}=\RR$, by definition.

\subsection{Determinant line of a direct sum}\Label{SS:directsum}
Given two objects $M$ and $N$ in $\c$, and a pair $\<\ ,\ \>_M$ and $\<\ ,\ \>_N$ of
$\c$-admissible scalar products on $M$ and $N$ correspondingly, we may obviously define
the scalar product $\<\ ,\ \>_M \oplus \<\ ,\ \>_N$ on the direct sum $M\oplus N$. This
defines an isomorphism
\begin{equation}\Label{E:directsum}
  \begin{CD}
        \phi: \det{M}\otimes\det{N} @>{\sim}>> \det(M\oplus N).
  \end{CD}
\end{equation}
Using \reft{FKdet} (properties (a) and (f)), it is easy to show that {\em this
homomorphism is canonical}, that is, it does not depend on the choice of the metrics $\<\
,\ \>_M$ and $\<\ ,\ \>_N$. From the description given above it is clear that the
homomorphism \refe{directsum} preserves the orientations.

\subsection{Push-forward of determinant lines}\Label{SS:push}
Note that, {\em any isomorphism $f:M\to N$ between objects $M,N$ of $\c$ induces
canonically an isomorphism of the determinant lines
\begin{equation}\Label{E:push}
    f_\ast:\det{M}\to\det{N}.
\end{equation}
Moreover, the induced map $f_\ast$ preserves the orientations of the determinant lines.}

Indeed, if $\<\ ,\ \>_M$ is an $\c$-admissible scalar product on $M$, then we set
\begin{equation}\Label{E:push2}
    f_\ast(\<\ ,\ \>_M) \ = \ \<\ ,\ \>_N,
\end{equation}
where $\<\ ,\ \>_N$ is the scalar product on $N$ given by
\begin{equation}\Label{E:pushsp}\notag
    \<v,w\>_N \ = \ \<f^{-1}(v),f^{-1}(w)\>_M
    \qquad\text{for}\quad v,w\in N
\end{equation}
(this scalar product is $\c$-admissible since $f\in \Hom_\c(M,N)$ is an isomorphism).
This definition does not depend on the choice of the scalar product $\<\ ,\ \>_M$ on $M$:
if we have a different $\c$-admissible scalar product $\<\ ,\ \>_M^\prime$ on $M$, where
$\<v,w\>_M^\prime = \<Av,w\>_M$ with $A\in\GL(M)$, then the induced scalar product on $N$
will be $\<v,w\>_N^\prime\ =\ \<({f}^{-1}Af)v,w\>_N$ and our statement follows from
property (a) of the Fuglede-Kadison determinant, cf. \reft{FKdet}.

\subsection{Calculation of the push-forward}\Label{SS:fMtoN}
Suppose $f:M\to N$ is an isomorphism and fix $\c$-admissible scalar products  $\<\ ,\
\>_M, \ \<\ ,\ \,\>_N$ on $M$ and $N$ respectively  (note that we don't assume any more
that the equality \refe{push2} is satisfied). Let $f^*:N\to M$ be the adjoint of $f$ with
respect to these scalar products. Then, for every $v,w\in N$ we have
\begin{equation}\Label{E:fMtoN1}\notag
    \<\, f^{-1}(v) ,f^{-1}(w)\, \>_M
    \ = \ \<\, (ff^*)^{-1}(v) ,w\, \>_N.
\end{equation}
Comparing this equation with \refe{push2} and \refe{equiv}, we obtain
\begin{equation}\Label{E:fMtoN}
    f_*\big(\, \<\ ,\ \>_M) \ = \ \sqrt{\Det_{\tau}(ff^\ast)}\cdot\<\ ,\ \>_N
\end{equation}
In particular, in the case $M= N$ we obtain the following

\prop{detf2}
If $f:M\to M$ is an automorphism of $M\in \Ob(\c)$, then the induced homomorphism
$f_\ast:\det{M}\to \det{M}$ coincides with the multiplication by $\Det_{\tau}(f)\in
\RR^{>0}$.
\eprop

\subsection{Composition of the push-forwards}\Label{SS:comppush}
It is obvious from the definition, that the construction of the push-forward is
{\em functorial}: if $f:M\to N$ and $g:N\to L$ are two isomorphisms between
objects of $\c$ then
\begin{equation}\Label{E:comppush}
    (g\circ f)_\ast\ =\  g_\ast\circ f_\ast.
\end{equation}

\prop{seqofmod}
Any exact sequence
\begin{equation}\Label{E:seqofmod}
  \begin{CD}
        0\to M^\prime@>{\alpha}>>M@>{\beta}>>M^{\prime\prime}\to 0
  \end{CD}
\end{equation}
of objects in $\c$ determines canonically an isomorphism
\begin{equation}\Label{E:seqofmod2}
    \psi_{\alp,\bet}:\, \det{M} \ \to  \ \det{M'}\otimes\det{M''}
\end{equation}
preserving the orientations of the lines. \eprop \prf The map
$\psi_{\alp,\bet}$ is constructed as follows. Any admissible scalar product
$\<\cdot , \cdot \>_M$ on $M$ defines admissible scalar products $\<\cdot ,
\cdot \>_{M'}$ on $M'$ and $\<\cdot , \cdot \>_{M''}$ on $M''$ as follows:
$$\<x',y'\>_{M'} = \<\alpha(x'), \alpha(y')\>_M, \quad
\<x'',y''\>_{M''} = \<\gamma(x''), \gamma(y'')\>_M$$ for $x',y'\in M'$ and $x'',y''\in
M''$. Here $\gamma: M''\to M$ is such that $\beta\circ \gamma =1_{M''}$ and
$\<\alpha(x'), \gamma(x'')\>_M =0$ for any $x'\in M'$ and $x''\in M''$. These conditions
determine $\gamma$ uniquely; in fact $\gamma$ identifies $M''$ with the orthogonal
complement $\im(\alpha)^\perp$ with respect to $\<\cdot, \cdot\>_M$.

Next we show that $\psi_{\alpha, \beta}$ is well-defined, i.e. the class of $\<\, \,
,\,\,\>_{M'}\oplus \<\, \, ,\,\,\>_{M''}$ in the determinant line $\det(M'\oplus M'')$
depends only on the class of the scalar product $\<\,\,,\,\,\>_M$ in $\det(M)$. Consider
an automorphism $h: M\to M$ and the induced scalar product $\<x,y\>'_{M}=\<h(x),h(y)\>_M$
on $M$. Let $\gamma: M''\to M$ be the splitting as above. The splitting $ \bar \gamma:
M''\to M$ corresponding to the the new scalar product $\<\cdot, \cdot\>'_{M}$ equals
$\bar \gamma =\gamma +\alpha\circ f$, where $f: M''\to M'$ is a morphism which is
uniquely determined by the data. One has
\begin{eqnarray}\label{stam1}
[\alpha, \bar \gamma] = [\alpha, \gamma]\circ \left[\begin{matrix} 1& f\\
0& 1
\end{matrix}
\right]. \end{eqnarray} Here $[\alpha, \gamma]$ denotes a morphism $M'\oplus
M''\to M$ whose restrictions to $M'$ and $M''$ are $\alpha$ and $\gamma$,
correspondingly. The symbol $[\alpha, \bar \gamma]$ has a similar meaning. The
matrix
$\left[\begin{matrix} 1& f\\
0& 1
\end{matrix}
\right]$ represents a morphism $M'\oplus M''\to M'\oplus M''$. Let $\xi, \xi'\in
\det(M'\oplus M'')$ be such that
$$[\alpha, \gamma]_\ast
(\xi) =\<\, \, ,\, \, \>_M,\quad [\alpha, \bar \gamma](\xi') = \<\, \, ,\, \, \>'_M.$$
Using the properties of the push-forwards applied to (\ref{stam1}) and \reft{FKdet}, we
find
$$[\alpha, \bar \gamma](\xi) = \<\, \, ,\, \, \>_M= \Det_\tau(h)^{1/2}\<\, \, ,\, \, \>'_M.$$
Thus, $\xi'=\Det_\tau(h)^{-1/2}\cdot \xi\, \in\,  \det(M'\oplus M'')$. In the determinant
line $\det(M)$ one has the relation $\<\, \, ,\, \, \>'_M=\Det_\tau(h)^{-1/2}\cdot \<\,
\, ,\, \, \>_M$; hence it follows that our definitions
 $\xi=\psi_{\alpha, \beta}(\<\, \, ,\, \, \>_M)$ and
$\xi'=\psi_{\alpha, \beta}(\<\, \, ,\, \, \>'_M)$ are compatible. \eprf

The following lemma will be used in \refs{detefc}. \lem{h} Under the conditions
of \refp{seqofmod} let $h:M\to M$ be an automorphism of $M$. Then the
automorphism $\psi_{h\circ\alp,\bet\circ{}h^{-1}}$ corresponding to the exact
sequence
\begin{equation}\Label{E:seqh}
    \begin{CD}
      0\to M^\prime@>{h\circ\alpha}>>M@>{\beta\circ h^{-1}}>>M^{\prime\prime}\to 0
    \end{CD}
\end{equation}
coincides with $\Det_\tau(h)^{-1}\cdot \psi_{\alp,\bet}$.
 \elem \prf
Fix an admissible scalar product $\<\cdot , \cdot \>_M$ on $M$. Let $\gamma: M''\to M$ be
such that $\beta\circ \gamma =1_{M''}$ and $\<\alpha(x'), \gamma(x'')\>_M =0$ for any
$x'\in M'$ and $x''\in M''$. Then
$$\psi_{\alpha, \beta}(\<\,\, ,\,\, \>_M)=\xi\in
\det(M'\oplus M'')$$ where $\xi$ satisfies
\begin{eqnarray}
[\alpha,\gamma]_\ast(\xi)=\<\,\, ,\,\, \>_M. \label{one}\end{eqnarray} This is
the definition of $\psi_{\alpha, \beta}$, see above. Now consider the sequence
(\ref{E:seqh}). Fix the following scalar product on $M$:
$$\<x,y\>'_M = \<h^{-1}(x), h^{-1}(y)\>_M, \quad x, y, M.$$
The corresponding splitting is $h\circ \gamma: M''\to M$. Hence
\begin{eqnarray}
\psi_{h\circ\alp,\bet\circ{}h^{-1}}(\<\,\, ,\,\, \>'_M)=\xi', \label{two}
\end{eqnarray}
where \begin{eqnarray} [h\circ \alpha, h\circ \gamma]_\ast(\xi') = \<\,\, ,\,\,
\>'_M.\label{three}\end{eqnarray} In the determinant line $\det(M)$ we have
\begin{eqnarray} \<\,\, ,\,\,\>'_M = \Det_\tau(h)\cdot \<\,\, ,\,\,
\>_M.\label{four}\end{eqnarray} Using the properties of the push-forwards we
obtain \begin{eqnarray} [h\circ \alpha, h\circ \gamma]_\ast(\xi') = h_\ast\circ
[\alpha,\gamma]_\ast(\xi') = \Det_\tau(h)\cdot
[\alpha,\gamma]_\ast(\xi').\label{five}\end{eqnarray} Combining (\ref{one}),
(\ref{two}), (\ref{three}), (\ref{four}), (\ref{five}) one obtains $\xi=\xi'$.
This proves that $\psi_{h\circ\alp,\bet\circ{}h^{-1}}= \Det_\tau(h)^{-1}\cdot
\psi_{\alp,\bet}$ as stated.

\eprf

\section{Abelian extension of a von Neumann category}\Label{S:extended}

{}M. Farber in \cite{Farber96} and \cite{Farber98} constructs an abelian
category called the {\em extended category of Hilbertian modules}. This abelian
category contains a given von Neumann category $\c$ as a full subcategory of
its projective objects. The construction of the extended category of
\cite{Farber96}, \cite{Farber98} was inspired by the earlier work of P. Freyd
\cite{Freyd} on embedding of additive categories into abelian categories.

Any object of the extended category splits naturally into projective and
torsion components; the projective part has the von Neumann dimension as an
invariant while the Novikov-Shubin invariant \cite{NovSh86} depends only on the
torsion part.

A brief description of Farber's construction is given below. We refer to the
paper \cite{Farber98} for more details.

\defe{extended}
Let $\c$ be a Hilbertian von Neumann category. The {\em extended category}
$\ec$ of $\c$ is defined as follows: An object of the category $\ec$ is defined
as a morphism $(\alpha:A^\prime\to A)$ in the original category $\c$. Given a
pair of objects $\calX=(\alpha: A^\prime \to A)$ and $\calY=(\beta:B^\prime\to
B)$ of $\ec$, a {\em morphism} $\calX\to\calY$ in category $\ec$ is an
equivalence class of morphisms $f:A\to B$ of category $\c$ such that
$f\circ\alpha= \beta\circ f'$ for some morphism $f':A^\prime \to B^\prime$ in
$\c$. Two morphisms $f_1:A\to B$ and $f_2:A\to B$ of $\c$ represent {\em
identical morphisms $\calX\to\calY$ of $\ec$} iff $f_1-f_2 = \beta\circ F$ for
some morphism $F:A\to B^\prime$ of category $\c$. The morphism $\calX\to\calY$,
represented by $f:A\to B$, is denoted by
\begin{equation}\Label{E:[f]}
    [f]:\, (\alpha:A^\prime\to A)\ \to\ (\beta:B^\prime\to B)
    \qquad \text{or by}\qquad [f]:\, \calX\to\calY.
\end{equation}
Composition of morphisms in $\ec$ is defined as composition of the
corresponding morphisms $f$ in the category $\c$. \edefe

The category $\ec$ is an abelian category, cf. \cite[Proposition~1.7]{Farber98}.

It is shown in \cite[Section~1.4]{Farber98}, that any object $\calX$ of $\ec$ is
isomorphic in $\ec$ to an object $(\alpha:A'\to A)$, where the morphism $\alpha$ is
injective.

\subsection{Projective and torsion objects}\Label{SS:proj-tor}
There is a full embedding of $\c$ in $\ec$ which takes an object $A$ of $\c$ to
the zero morphism $(0\to A)$ of $\ec$ and a $\c$-morphisms $f:A\to B$ to the
morphism $[f]:(0\to A)\to(0\to B)$ of $\ec$. An object of $\ec$ is projective
if and only if it is isomorphic to a Hilbertian module, i.e. to an object of
$\c$, cf. \cite[Proposition~1.9]{Farber98}.

Given an object $\calX=(\alp:A'\to A)$ in the extended category $\ec$, its {\em
projective part} is defined as the following object of $\c$
\begin{equation}\Label{E:projpart}
    \calP(\calX) \ := \ A/\cl(\im(\alp)).
\end{equation}
Clearly, any morphism $[f]:\calX\to \calY$ of $\ec$ induces a morphism
$\calP(f):\calP(\calX)\to \calP(\calY)$ between the projective parts. Thus we
have a well defined functor
\begin{equation}\Label{E:functorproj}
    \calP:\, \ec \ \to \ \c.
\end{equation}

An object $\calX=(\alpha:A'\to A)$ of the extended category $\ec$ is called
{\em torsion} iff the image of $\alpha$ is dense in $A$. One denotes by $\tca$
the full subcategory of $\ec$ generated by all torsion objects. $\tca$ is
called {\em the torsion subcategory of $\ec$}. It is shown in
\cite[\S3]{Farber98}, that: {\em if $\c$ is a finite von Neumann category, then
$\tca$ is an abelian category.}

Objects of $\tca$ are called {\em torsion Hilbertian modules}.

There is a functor
\begin{equation}\Label{E:functortor}
    \calT:\, \ec \ \to \ \tca, \qquad
    \calT:\, (\alp:A'\to A) \ \mapsto \ (\alp:A'\to \cl(\im(\alp))).
\end{equation}
Given $\calX \in \ob(\ec)$, the object $\calT(\calX)$ is called {\em the
torsion part of $\calX$}.

It is shown in \cite[\S3.4]{Farber98} that {\em every object $\calX\in \ec$ is
a direct sum of its torsion and projective parts},
\begin{equation}\Label{E:tor+proj}
    \calX \ = \ \calT(\calX)\oplus \calP(\calX).
\end{equation}

\subsection{Isomorphisms in the extended category}\Label{SS:isomorphism}
Let $\calX= (\alp:A'\to A)$ and $\calY= (\bet:B'\to B)$ be two objects of the
extended category $\ec$. From \refd{extended} we see that a morphism $f\in
\Hom_\c(A,B)$ induces an isomorphism in $\ec$ if and only if there exist
$\c$-morphisms $\bet, f, f', g, g', F, G$ such that the following diagrams
commute
\begin{equation}\Label{E:fgenerated}
    {\begin{CD}
    A' & @>\alp>>& A\\
    @Vf'VV&   &  @VV{f}V\\
    B' & @>>\bet> & B
    \end{CD}}
    \qquad\quad
    \begin{CD}
    A' & @>\alp>>& A\\
    @Ag'AA&   &  @AAgA\\
    B' & @>>\bet> & B
    \end{CD}
    \qquad\quad
    \begin{CD}
    A' & @>\alp>>&& A&&\\
      & &{}^{{}^F}\text{\Large $\swarrow$}   &      && @VV{\id_A-gf}V\\
    A' & @>>\alp> && A&&
    \end{CD}
    \qquad\qquad
     \begin{CD}
    B' & @>\bet>>&& B&&\\
      & &{}^{{}^G}{\text{\Large$\swarrow$}}   &  &&@VV{\id_B-fg}V\\
    B' & @>>\bet> && B&&
    \end{CD} \qquad{}
\end{equation}

\lem{sequence} Let $\calX=(\alp:A'\to A)$ and $\calY=(\bet:B'\to B)$ be two
objects of $\ec$ which are isomorphic is $\ec$ and such that the morphisms
$\alp$ and $\bet$ are injective. Let $f:A\to B$ be a morphism in $\c$ which
induces an isomorphism $[f]:\calX\to \calY$. Let $f':A'\to B'$ be as in the
left diagram of \refe{fgenerated}. Then the following sequence of morphisms of
$\c$ is exact\footnote{Since the morphisms $\alp$ and $\bet$ are injective the
morphism $f'$ and the sequence \refe{sequence} are completely determined by
$f$.}:
\begin{equation}\Label{E:sequence}
    \begin{CD}
        0\to \ A' \ @>{(f',\, \alp)}>> \ B'\oplus A @>{\bet-f}>> \
          B \ \to \ 0.
    \end{CD}
\end{equation}
\elem \prf We shall apply Proposition 1.6 from \cite{Farber98} which describes
explicitly the kernels and cokernels in $\ec$. Clearly, $[f]:\calX \to \calY$
being an isomorphism is equivalent to the vanishing of its kernel and cokernel
in $\ec$. The cokernel of $[f]$ in $\ec$ equals $((\beta, -f): B'\oplus A\to
B)$ and its vanishing means that the $\calC$-morphism $(\beta, -f): B'\oplus
A\to B$ is onto; this is equivalent to the exactness of (\ref{E:sequence}) at
the last term. The kernel of $[f]$ in $\ec$ equals $((f',\alpha): A'\to P)$
where $P$ denotes the kernel of $\beta+f: B'\oplus A \to B$ (here we use the
assumption that $\beta$ is injective). Hence vanishing of the kernel of $[f]$
is equivalent to the exactness of (\ref{E:sequence}) in the middle term. The
exactness of (\ref{E:sequence}) in the first term follows from the injectivity
of $\alpha$.
 \eprf
\begin{remark}{remarknew}
One may think about an object $\calX=(\alpha: A'\to A)$ of $\ec$ as of a short
chain complex in $\calC$ where $\alpha$ is the differential. If $\alpha$ is
injective then this chain complex has cohomology in a single dimension (i.e. it
is a resolution). A morphism in $\ec$
\begin{eqnarray}\label{E:morphism}
{\begin{CD}
    A' & @>\alp>>& A\\
    @Vf'VV&   &  @VV{f}V\\
    B' & @>>\bet> & B
    \end{CD}}
    \end{eqnarray}
 is then a morphism of chain complexes. The chain complex
(\ref{E:sequence}) is the cone of this chain map. Hence the lemma above can be
rephrased by saying that a morphism of the extended category is an isomorphism
if and only if its cone is acyclic.
\end{remark}


\section{Determinant lines of objects in the extended category}\Label{S:detefc}

Throughout this section $\c$ denotes a Hilbertian von Neumann category endowed
with a finite trace, cf. \refss{traces}. We denote by $\ec$ the extended
category of $\c$, cf. \refs{extended}.

\defe{detmorphism}
Let $\calX=(\alp:A'\to A)\in \ob(\ec)$. The determinant line of $\calX$ is
\begin{equation}\Label{E:detecf}
    \det\calX \ :=
    \det{A}\otimes (\det{A'/\ker(\alpha)})^*.
\end{equation}
\edefe

An element of the determinant line $\det \calX$ can be represented by a symbol
of the form
\begin{eqnarray}\label{E:symbol}
\lambda\cdot \frac{\langle\, ,\, \rangle}{\langle\, ,\, \rangle^\prime},
\end{eqnarray} where $\lambda\in \mathbb R$ and $\langle\, ,\, \rangle$,
$\langle\, ,\, \rangle^\prime$ are $\calC$-admissible scalar products on $A$,
and on $A'/\ker(\alpha)$ correspondingly. The equivalence relation between
symbols (\ref{E:symbol}) follows from (\ref{E:equiv}) and from the usual rules
of fractions.

If $\calX=(0\to A)$ is a projective object in $\ec$, cf. \refss{proj-tor}, then
the determinant line $\det\calX$ coincides with $\det{A}$.

Let $\calX=(\alp:A'\to A)$ and $\calY=(\bet:B'\to B)$ be two objects of the
extended category $\ec$ such that the morphisms $\alp$ and $\bet$ are
injective. Consider a morphism $[f]:\calX\to \calY$ represented by diagram
(\ref{E:morphism}). Assuming that $[f]$ is an isomorphism in $\ec$ we define
the push-forward map $[f]_*:\det\calX\to \det\calY$ associated with $[f]$ as
follows. Fix arbitrary $\c$-admissible scalar products $\br_A$, $\br_{A'}$,
$\br_{B'}$ on $A$,$A'$ and $B'$ correspondingly. Consider the exact sequence
(\ref{E:sequence}). The orthogonal complement to the image of $(f',\alpha)$
(with respect to the scalar product $\langle \, , \, \rangle_{B'}\oplus \langle
\, , \, \rangle_A$) is mapped by $\beta -f$ isomorphically onto $B$. This
isomorphism together with the restriction of the scalar product $\br_{B'}\oplus
\br_A$ onto ${\im(f',\alpha)^\perp}$ determine a $\calC$-admissible scalar
product $\langle \, , \, \rangle_B$ on $B$. We set
$$[f]_\ast\left(\frac{\langle \, , \, \rangle_A}{\langle \, , \, \rangle_{A'}}\right)=
\frac{\langle \, , \, \rangle_B}{\langle \, , \, \rangle_{B'}}\, \in\,
\det\calY.$$ Recall that the symbol $\displaystyle{\frac{\br_A}{\br_{A'}}}$
denotes a nonzero element of the line $\det \calX$.

 \lem{[f]} The push-forward map $[f]_*:\det\calX\to
\det\calY$ is well-defined. It depends only on the class $[f]$ in
$\Hom_{\ec}(\calX,\calY)$. \elem

\prf Assume that $\tilf:A\to B$ is another morphism such that $[\tilf]= [f]$.
Then there exists a morphism $R:A\to B'$ such that $\bet\circ{}R= f-\tilf$. Let
$\tilf':A'\to B'$ satisfy $\bet\circ\tilf'= \tilf\circ\alp$. Then
\[
    \bet\circ(f'-\tilf') \ = \ (f-\tilf)\circ\alp \ = \ \bet\circ R\circ\alp,
\]
and, since $\bet$ is injective, $f'-\tilf'= R\circ\alp$.

Set
\[
    h \ = \ \begin{pmatrix}
                \id_{B'} & R\\
                0   & \id_{A}
            \end{pmatrix}:\, B'\oplus A \ \to \  B'\oplus A.
\]
Then
 \(
    h^{-1}=\left(
    \begin{smallmatrix}
        \id_{B'} & -R\\
        0       & \id_{A}
    \end{smallmatrix}\right):B'\oplus A \to B'\oplus A
 \)
and
\begin{align}
    (\tilf',-\alp) \ &= \ (f'-R\circ\alp,-\alp) \ = \ h\circ(\tilf',-\alp),\notag\\
    \qquad
    \bet + \tilf \ &= \ \bet+f -\bet\circ R \ = \ (\bet+\tilf)\circ h^{-1}.\notag
\end{align}
Since $\Det_{\tau}(h)= \Det_{\tau}(h^*)= 1$ by \refe{detM+N} our claim follows
now from \refl{h}. \eprf

To define the push-forward map $[f]_\ast: \det \calX \to \det\calY$ without
assuming that $\alpha$ and $\beta$ are injective one observes that the
determinant lines of $\calX$ and of $\calX'=(\alpha: A'/\ker \alpha\to A)$ are
identical and similarly the determinant line of $\calY$ and of $\calY'=(\beta:
B'/\ker\beta\to B)$ are identical. Any isomorphism $[f]:\calX\to \calY$ given
by the diagram (\ref{E:morphism}) determines an isomorphism $[f]: \calX'\to
\calY'$ given by
\begin{eqnarray}
{\begin{CD}
    A'/\ker \alpha & @>\alp>>& A\\
    @Vf'VV&   &  @VV{f}V\\
    B'/\ker \beta & @>>\bet> & B.
    \end{CD}}
    \end{eqnarray}

The next claim follows directly from the definitions:

\prop{functoriality}
Let $[f]:\calX\to \calY$ and $[g]:\calY\to \calZ$ be isomorphisms. Then
\begin{equation}\Label{E:functoriality}\notag
    [g]_*\circ [f]_* \ = \ \big(\, [g]\circ[f]\, \big)_*.
\end{equation}
\eprop

\subsection{Determinant line of a direct sum}\Label{SS:propdet}
Given two objects $\calX, \calY\in \ec$ there is a canonical isomorphism
\begin{equation}\Label{E:directsum2}
    \phi: \det\calX\otimes\det\calY \to    \det(\calX\oplus \calY),
\end{equation}
compare (\ref{E:directsum}).

In particular one obtains that $\det\calX$ is canonically isomorphic to the
product of the determinants of the projective and the torsion parts:
$
    \det\calX \ \simeq \ \det\calP(\calX)\otimes \det\calT(\calX).
$

\lem{kerXcoker}\label{E:kerX} For any morphism $[f]:\calX\to \calY$ of $\ec$
there is a canonical isomorphism
\begin{equation}\Label{E:kerXcoker}\notag
  \begin{CD}
    \phi_f:\,\det\calY\otimes (\det\calX)^* \ @>{\sim}>> \
      \det\coker([f])\otimes(\det\ker([f]))^*.
  \end{CD}
\end{equation}
\elem
\prf
Let $\calX=(\alp:A'\to A)$ and $\calY=(\bet:B'\to B)$. Without loss of generality we can
assume that $\alp$ and $\bet$ are injective.

Let $f:A\to B$ be a representative of $[f]$. Then there is a unique $f':A'\to B'$ such
that $f\circ\alp= \bet\circ{}f'$. Recall that the kernel and the cokernel of $[f]$ are
described in Proposition~1.6 of \cite{Farber98} as follows.
\begin{equation}\Label{E:kercoker}
    \coker([f]) \ = \ ((\bet, -f):B'\oplus A\to B),
    \qquad
    \ker([f]) \ = \ (i:P'\to P),
\end{equation}
where
\begin{equation}\Label{E:P'-P}
    P \ = \ \ker((\bet, -f):\, B'\oplus A\to B),
    \qquad
    P' \ = \ \ker((\bet, \, -f\circ \alpha):\, B'\oplus A'\to B),
\end{equation}
and $i:P'\to P$ is the restriction of the map ${\id}\oplus\alp$ to $P'$. The
assumption that $\bet$ is injective allows to simplify the description of $P'$.
Indeed, by \refe{P'-P}, $P'$ consists of those pairs $(a',b')\in A'\oplus{}B'$
such that $0 \ = \ f(\alp(a')) - \bet(b') \ = \ \bet(f'(a')-b'))$; the latter
is equivalent to $f'(a') = b'$. Thus there is a natural isomorphism
$A'\overset{\sim}{\to} P'$ given by $a'\mapsto (f'(b'), a')$, and the kernel of
$[f]$ can be rewritten as
$
    \ker([f]) \ = \  (\iot:A'\to P),
    \, \,
    \iot:a'\mapsto (f'(a'), \alp(a')).
$
Since $\alp$ is injective so is $\iot$. Thus we obtain
\begin{equation}\Label{E:detker}
    \det\ker([f])
    \ = \ \det{P}\otimes(\det{A'})^* \ = \
    \det(\ker(\bet\oplus-f))\otimes(\det{A'})^*.
\end{equation}

{}From \refe{kercoker} we obtain
\begin{multline}\Label{E:detcoker}
    \det\coker([f]) \ = \
    \det{B}\otimes(\det(B' \oplus A/\ker(\bet\oplus -f)))^* \ =
    \\ = \
    \det{B}\otimes(\det{B'})^*\otimes(\det{A})^*
        \otimes\det(\ker(f\oplus-\bet)).
\end{multline}

Combining \refe{kercoker}, \refe{detker} and \refe{detcoker} we obtain a
natural isomorphism
\begin{equation}\Label{E:kercoker2}
    \det\coker([f])\otimes(\det\ker([f]))^* \ \cong \
    \det{B}\otimes(\det{B'})^*\otimes(\det{A})^*
    \otimes\det{A'}.
\end{equation}
Since the RHS of \refe{kercoker2} is exactly $\det\calY\otimes(\det\calX)^*$
the lemma follows. \eprf

\section{$\tau$-trivial torsion objects and their determinant
lines}\label{triviality}

Let $\calC$ be a Hilbertian von Neumann category endowed with a finite trace
$\tau$. We have shown above how one may associate determinant lines to the
objects of the extended category $\ec$.

The trivial object of $\ec$ is represented by any morphism $\calX=(\alpha:A'\to
A)$ such that $\alpha$ is onto. Let us show that the determinant line of the
trivial object is {\it trivial} in the sense that it contains a canonical
nonzero element. The choice of such element allows to establish a unique
isomorphism to standard line $\mathbb R$ sending the canonical element to 1.
Let $\br$ be a $\calC$-admissible scalar product on $A$. The morphism $\alpha:
A'\to A$ determines an isomorphism $\alpha': A'/\ker \alpha \to A$. We set
$\langle x,y\rangle'=\langle \alpha'(x), \alpha'(y)\rangle$ where $x, y\in
A'/\ker \alpha$. The bracket $\br'$ is a $\c$-admissible scalar product on
$A'$. It is easy to see that the following nonzero element
$$\frac{\br}{\br'}\, \in \, \det\calX$$
is {\it canonical}, i.e. it does not depend on the choice of the initial scalar
product $\br$.

The purpose of this section is to show that such trivialization of the
determinant lines happens for a class of torsion objects which we call
$\tau$-{\it trivial}; the latter class depends on the choice of the trace
$\tau$.

Let $\calX = (\alpha:A'\to A)$ be a torsion object of $\ec$. This means that
the image of $\alpha$ is dense in $A$. Without loss of generality we may assume
that $\alpha$ is injective.

\defe{tauiso}
A morphism $\alpha: A'\to A$ is called a $\tau$-isomorphism if it is injective,
has dense image, and the integral
\begin{eqnarray}\label{integrall}
\int\limits_0^\infty \ln \lambda d\phi(\lambda) \, >\, -\infty
\end{eqnarray}
is finite.\edefe

Here $\phi(\lambda)$ denotes the spectral density function
$\phi(\lambda)=\tau(E_\lambda)$ of the self-adjoint operator
\begin{eqnarray}
(\alpha^\ast\alpha)^{1/2} = \int\limits_0^\infty \lambda dE_\lambda.
\label{specdec}\end{eqnarray} The adjoint $\alpha^\ast$ is calculated with
respect to a choice of $\calC$-admissible scalar products $\br$ on $A$ and
$\br_1$ on $A'$. The finiteness of the integral (\ref{integrall}) is
independent of this choice (see Proof of Proposition 3.2 in
\cite{CarFarMat97}). Note that the integral (\ref{integrall}) may only diverge
at point $\lambda=0$.

We shall use the following notation \begin{eqnarray}
\Det_\tau(\alpha^\ast\alpha) = \exp\left[2\int\limits_0^\infty
\ln(\lambda)d\phi(\lambda)\right]\label{FK}\end{eqnarray} where $\alpha$ and
$\phi(\lambda)$ are as above. It is an extension of the Fuglede - Kadison
determinant, see \cite{FK}, Lemma 5; this notation is also compatible with the
formula (\ref{exdet}) above.

\defe{tautriv} A torsion object $\calX$ is called $\tau$-trivial
if it is isomorphic to $(\alpha:A'\to A)$ where $\alpha$ is a
$\tau$-isomorphism. \edefe

\begin{lemma}{tautriv1}
The determinant line $\det \calX$ of any $\tau$-trivial torsion object $\calX$
is trivialized in the above sense, i.e. it contains a canonical nonzero
element.
\end{lemma}

\begin{proof} Let $\calX=(\alpha:A'\to A)$ with $\alpha$ being a $\tau$-isomorphism.
Fix $\calC$-admissible scalar products $\br$ on $A$ and $\br'$ on $A'$.
Consider the spectral decomposition (\ref{specdec}) where the adjoint
$\alpha^\ast$ is calculated with respect to the chosen scalar products. Let us
show that the following nonzero element
\begin{eqnarray}\label{element}
\frac{\br}{\br'}\cdot\sqrt{\Det_\tau(\alpha^\ast\alpha)}\, \in \, \det\calX
\end{eqnarray}
is independent of the choices of the scalar products $\br$ and $\br'$. Let
$\beta:A\to A$ and $\gamma: A'\to A'$ be positive self-adjoint invertible
$\calC$-morphisms. Consider the new scalar products $\br_1=\langle \beta\cdot,
\cdot\rangle$ and $\br_2=\langle \gamma\cdot, \cdot\rangle'$ on $A$ and on $A'$
correspondingly. The adjoint of $\alpha$ with respect to the new pair of scalar
products is $\gamma^{-1}\alpha^\ast\beta$. Using the multiplicativity property
of the Fuglede - Kadison determinant we find
\begin{eqnarray}\nonumber\Det_\tau(\gamma^{-1}\alpha^\ast \beta\alpha)^{1/2} =
\Det_\tau(\gamma)^{-1/2} \Det_\tau(\alpha^\ast\beta\alpha)^{1/2}=\\
=\,\Det_\tau(\gamma)^{-1/2} \Det_\tau(\alpha\alpha^\ast\beta)^{1/2}=
\Det_\tau(\gamma)^{-1/2}
\Det_\tau(\alpha\alpha^\ast)^{1/2}\Det_\tau(\beta)^{1/2}=\nonumber\\
=\,
\Det_\tau(\gamma)^{-1/2}\Det_\tau(\alpha^\ast\alpha)^{1/2}\Det_\tau(\beta)^{1/2}
\nonumber\end{eqnarray} Relation (\ref{E:equiv}) gives
$\br_1=\Det_\tau(\beta)^{-1/2}\cdot \br$ and $\br_2=
\Det_\tau(\gamma)^{-1/2}\cdot \br'$. Hence combining the above relations we
find that the new density
\begin{eqnarray} \frac{\br_1}{\br_2}\cdot
\Det_\tau(\gamma^{-1}\alpha^\ast \beta\alpha)^{1/2}\label{newelement}
\end{eqnarray}
equals (\ref{element}).
\end{proof}

\section{Determinant line of a chain complex}\Label{S:detcomp}

\subsection{}\Label{SS:graded}
Let $\calM= \oplus{}\calM^i$ be a graded object of $\ec$. The determinant line
of $\calM$ is defined as
\begin{equation}\Label{E:detgrobj}
    \det\calM \ := \ \otimes\, (\det\calM^i)^{(-1)^i}
\end{equation}
where $(\det\calM^i)^{-1}$ denotes the dual line to $\det\calM^i$.

Suppose now that
\begin{equation}\Label{E:complex}
  \begin{CD}
    (\calM,\p): \ \
    0 \ \to \calM^0 \ @>{\p_1}>> \calM^1 @>{\p_2}>> \cdots @>{\p_n}>> \calM^n \to 0
  \end{CD}
\end{equation}
is a chain complex in the abelian category $\ec$. We denote by
\begin{eqnarray}\label{extended}
    \calH^\ast(\calM) \ =\ \bigoplus_{i=0}^n\calH^i(\calM)
\end{eqnarray}
its cohomology. Even if the initial complex $\calM$ consisted of projective
objects (i.e. if $\calM^i\in \ob(\calC)$ for all $i$), the cohomology
$\H^\ast(\calM)$ may have nontrivial torsion, thus it is well-defined only as
an object of the extended category $\ec$, i.e. $\H^i(\calM)\in \ob(\ec)$. It is
called the {\it extended cohomology} of $\calM$. The notion of extended
cohomology was first introduced in \cite{Farber95}, \cite{Farber96}.

In the special case when the complex $\calM$ is projective, the extended
cohomology can explicitly be expressed by \begin{eqnarray} \H^i(\calM) =
(\partial_i:\, \calM^{i-1}\to \ker(\partial_{i+1})),
\label{extended1}\end{eqnarray} see \cite{Farber95}, \cite{Farber96}.

The projective part of the extended cohomology is then isomorphic to the {\em{}
reduced cohomology}
\begin{equation}\Label{E:reduced}
    H^i(\calM) \ := \ \ker(\p_{i+1})/\cl(\im(\p_{i}))\, \in\, \ob(\calC), \qquad i=0\nek n,
\end{equation}
of $(\calM,\p)$. The notion of reduced $L^2$-cohomology was originally
introduced by M.Atiyah \cite{Atiyah76}. The torsion part of the extended
cohomology is responsible for, so called, ``zero in the continuous spectrum
phenomenon". In particular, it completely describes the Novikov-Shubin
invariants of $(\calM,\p)$, cf. \cite[\S3]{Farber98}. The torsion part of
$\H^i(\calM)$ equals
$$\calT(H^i(\calM))\, =\, (\p_i:\,
\calM^{i-1}\to \cl(\im(\p_i)))$$ see \cite{Farber96}, \cite{Farber98}.

\prop{detofcoh} The chain complex \refe{complex} defines a canonical
isomorphism
\begin{equation}\Label{E:detofcoh}
    \nu_\calM:\,\det\calM \ \longrightarrow \ \det\calH^\ast(\calM),
\end{equation}
preserving the orientations of the lines. \eprop 
\prf We use here the notions \lq\lq image\rq\rq\, and \lq\lq kernel\rq\rq\ in
their categorical sense; we apply them to morphisms of the abelian category
$\ec$.

Let $Z^i(\calM)=\ker[\partial: \calM^i\to \calM^{i+1}]$ and
$B^i(\calM)=\im[\partial: \calM^{i-1}\to \calM^i]$ be the \lq\lq cycles\rq\rq \
and the boundaries, correspondingly. One obtains two familiar exact sequences
in $\ec$:
$$0\to B^i(\calM)\to Z^i(\calM) \to \H^i(\calM)\to 0$$
and
$$0\to Z^{i-1}(\calM)\to \calM^{i-1}\to B^i(\calM)\to 0.$$
The isomorphisms of \refp{seqofmod} applied twice (alternatively one may appeal
here to \refl{kerXcoker}), determines a canonical isomorphism
$$\det\H^i(\calM)\simeq (\det \calM^{i-1})^\ast\otimes \det Z^{i-1}(\calM)\otimes
\det Z^i(\calM).$$

Now one obtains
\begin{eqnarray*}
\det \calH^\ast(\calM) \, =\,  \prod_{i=0}^n (\det\H^i(\calM))^{(-1)^i} =
\\ =\, \prod_{i=0}^n(\det \calM^i)^{(-1)^i}\otimes \prod_{i=0}^n\left[\det
Z^i(\calM)\otimes \det Z^{i-1}(\calM)\right]^{(-1)^i} =\\
=\, \det(\calM) \otimes\prod_{i=0}^n\left[\det Z^i(\calM)\otimes \det
Z^{i-1}(\calM)\right]^{(-1)^i}.
\end{eqnarray*}
The main point of the proof is the observation that the product of the square
brackets in the last formula is a trivial line in the sense that it contains a
canonical nonzero element (denoted 1) and so it is canonically isomorphic to
$\mathbb R$. We see that the line $\det \H(\calM)$ is obtained from the line
$\det \calM$ by tensoring with a trivialized line containing a canonical
element which we shall denote by 1. The canonical isomorphism
(\ref{E:detofcoh}) is defined by the formula
$$\nu_\calM(x) =x\otimes 1.$$

 \eprf


\subsection{Direct sum of complexes}\Label{SS:tordsum}
Let $(\calM,\p^\calM)$ and $(\calN,\p^\calN)$ be two complexes of objects in
$\E(\c)$. Consider their direct sum
$(\calM\oplus\calN,\p^\calM\oplus\p^\calN)$. It follows from \refe{directsum2}
that
\begin{equation}\Label{E:directsum3}
    \nu_{\calM\oplus\calN} \ = \ \nu_\calM\otimes\nu_\calN.
\end{equation}

\subsection{Torsion of a complex of projective objects}\Label{SS:projcomp}
Suppose
\begin{equation}\Label{E:prcomplex}
  \begin{CD}
    (C,\p): \ \
    0 \ \to C^0 \ @>{\p_1}>> C^1 @>{\p_2}>> \cdots @>{\p_n}>> C^n \to 0
  \end{CD}
\end{equation}
is a complex in $\c$. We shall consider it as a complex of projective objects
of $\ec$. Fix $\calC$-admissible scalar products $\<\ ,\ \>_i$ on $C^i$. They
determine elements $\sig_i\in \det{C^i}$ and
$$\sigma=\prod_{i=0}^n \sigma_i^{{(-1)}^i}\, \in \det C$$

\defe{torsion}
The positive element
\begin{equation}\Label{E:torsion}
    \rho_C \ := \ \nu_C(\sig) \ \in \det\calH^\ast(C)
\end{equation}
is called the torsion of the complex $C$. \edefe 
\rem{torsion} The torsion $\rho_C$ depends on the choice of the scalar products
$\<\,,\,\>_i$. In \cite{CarFarMat97}, the torsion was defined more invariantly
as the element of the line $(\det{C})^*\otimes\det\H^\ast(C)$ canonically
determined by the map $\nu_C$. In this paper we adopt a less invariant
definition \refe{torsion}, since it is more consistent with the usual notions
of combinatorial and de Rham torsions. Note that, using our notation, the
torsion of \cite{CarFarMat97} can be written as $
    \sig^*\otimes \rho_C,
$
where $\sig^*$ is the unique element of $(\det{C})^*$ such that $\<\sig^*,\sig\>=1$.
\erem

\subsection{Torsion of an acyclic complex of projective objects}\Label{SS:acprojcomp}
Assume now that the complex \refe{prcomplex} is acyclic (by this we mean that
$\im\p_{i-1}= \ker\p_i$. The acyclicity implies in particular that $\im(\p_i)$
is a closed subspace of $C^i$). We now calculate the torsion $\rho_C=
\nu_C(\sig)\in \det\calH(C)\simeq \RR$ which is in this situation a positive
real number.

Let $\p_i^*: C^i\to C^{i-1}$ denote the adjoint of $\p_i$ with respect to the chosen
scalar products and let
\[
    \Del_i  \ := \ \p_{i+1}^*\p_{i+1}+\p_{i}\p_{i}^*:\, C^i \to C^i
\]
be the Laplacian. Since the complex $(C,\p)$ is acyclic, $\Del_i$ is a positive
invertible operator. We denote by $\Det_{\tau}(\Del_i)$ the Fuglede-Kadison determinant
of $\Del_i$, cf. \refss{Fkdet}.

\prop{projcomp} Assume that the chain complex \refe{prcomplex} is acyclic. Then
its torsion is given by the formula
\begin{equation}\Label{E:projcomp}
    \rho_C \ = \
    \exp\left(\frac12\, \sum_{i=0}^n\, (-1)^ii\log\Det_{\tau}(\Del_i)\right).
\end{equation}
\eprop
\prf
Set
\[
    A^i \ := \ \ker(\p_{i+1}) \ = \ \im(\p_i);
    \qquad
    B^i \ := \ \ker(\p_i^*) \ = \ \im(\p_{i+1}^*)
        \ = \ \big(\ker(\p_{i+1})\big)^{\perp}.
\]
Then $C^i=A^i\oplus B^i$ and the complex \refe{prcomplex} decomposes into
direct sum of the following \lq\lq short\rq\rq \ chain complexes
\[
 \begin{CD}
    (\calM_i,\p):\ \
      0 \ \to \ \cdots 0 \ \to \ B^{i-1} @>{\p_i|_{B^{i-1}}}>>
      \ A^i \ \to  \ 0 \ \to \cdots\to \ 0,
 \end{CD}
\]
where $i=1, \nek, n$. The scalar products $\<\ ,\ \>_i$ on $C^i$ $(i=0\nek n)$
induce scalar products on $A^i$ and $B^i$ and hence the elements $\mu_i\in
\det{\calM_i}$. Clearly,
\begin{equation}\Label{E:tor=prod}
    \sig \ =  \ \mu_1\otimes\cdots\otimes \mu_n,
    \qquad
    \nu_C \ = \ \nu_{\calM_1}\otimes \cdots\otimes \nu_{\calM_n}.
\end{equation}

{}From \refe{fMtoN} we see that
\begin{equation}\Label{E:torCi}
    \rho_{\calM_i} \ = \ \nu_{\calM_i}(\mu_i)
    \ = \ \left(\Det_{\tau}(\p_i^*\p_i|_{B^{i-1}})\right)^{(-1)^i/2}.
\end{equation}
Combining \refe{tor=prod} and \refe{torCi} with \refe{torsion} we obtain
\begin{multline}\Label{E:rho-sig}
    \rho_C \ = \
     \Det_{\tau}(\p_1^*\p_1|_{B^0})^{-1/2}\cdot\Det_{\tau}(\p_2^*\p_2|_{B^1})^{1/2}
    \cdots \Det_{\tau}(\p_{n}^*\p_{n}|_{B^{n-1}})^{(-1)^n/2}
    \\ = \
    \exp\Big(\frac12\, \sum_{i=1}^n\, (-1)^i
        \log\Det_{\tau}\big(\p_i^*\p_i|_{B^{i-1}}\, \big)\Big).
\end{multline}

It remains to show that the right hand sides of \refe{projcomp} and
\refe{rho-sig} coincide. The Laplacian has the following matrix form with
respect to the decomposition $C^i=A^i\oplus{}B^i$:
\[
    \Del_i \ = \ \begin{pmatrix}
                    \p_i\p_i^*&0\\
                    0&\p_{i+1}^*\p_{i+1}
                  \end{pmatrix}.
\]
Hence, using \refe{detM+N}, we obtain
\begin{equation}\Label{E:detAB}
    \Det_{\tau}(\Del_i) \ = \
    \Det_{\tau}(\p_i\p_i^*|_{A^i})\cdot\Det_{\tau}(\p_{i+1}^*\p_{i+1}|_{B^i}).
\end{equation}
The restriction of $\p_i$ to $B^{i-1}$ maps it isomorphically onto $A^{i}$ and we have
\begin{equation}\Label{E:pp*-p*p}\notag
    \p_i|_{B^{i-1}} \, \big( \p_i^*\p_i|_{B^{i-1}}\big) \ = \
    \big(\p_i\p_i^*|_{A^{i}}\big)\, \p_i|_{B^{i-1}}.
\end{equation}
Hence,
\begin{equation}\Label{E:Det=Det}
    \Det_{\tau}(\p_{i}^*\p_{i}|_{B^{i-1}}) \ = \ \Det_{\tau}(\p_i\p_i^*|_{A^i})
\end{equation}
by \reft{FKdet}(a). Using \refe{detAB} and \refe{Det=Det} we get
\begin{multline}\Label{E:DetA}
    \sum_{i=0}^n\, (-1)^ii\log\Det_{\tau}\Del_i \ = \
    \sum_{i=0}^{n-1}\, (-1)^ii\, \big(\, \log\Det_{\tau}(\p_{i}^*\p_{i}|_{B^{i-1}})
            + \log\Det_{\tau}(\p_{i+1}^*\p_{i+1}|_{B^{i}}) \, \big)
    \\ = \
    \sum_{i=1}^{n}\, (-1)^i\,\log\Det_{\tau}(\p_{i}^*\p_{i}|_{B^{i-1}}).
\end{multline}
Combining \refe{rho-sig} and \refe{DetA} we obtain \refe{projcomp}.
\eprf

\section{Torsion and exact sequences}\Label{S:exseq}

In this section we study the relationship between the torsions of the terms of
a short exact sequence of complexes. Roughly speaking the main result states
that the torsion of the middle term equals to the product of the torsions of
the two other terms.  This result will be used in \refs{combtor} to show that
the combinatorial torsion of a flat Hilbertian bundle is invariant under
subdivisions. It will be also used in \refs{chmu} to prove that the ratio of
the de Rham and the combinatorial torsions equals the relative torsion
introduced in \cite{CarMatMis99}.

\subsection{The setting}\Label{SS:setexceq}
Let $\c$ be a Hilbertian von Neumann category endowed with a finite trace
$\tau$.  Let
\begin{equation}\Label{E:exseq}
 \begin{CD}
    0 \ \longrightarrow \ L @>{\alp}>> M @>\bet>> N \ \longrightarrow 0.
 \end{CD}
\end{equation}
be an exact sequence of chain complexes $(L,\p^L)$, \ $(M,\p^M)$ and $(N,\p^N)$
of objects of $\c$. Let
\[
    \psi: \det M \ \overset{\sim}\longrightarrow \
    \det L\otimes \det N
\]
be the isomorphism determined by \refe{exseq}, cf. \refp{seqofmod}.

Choose an identification of $M$, as a graded Hilbertian module, with the direct
sum of $L$ and $N$. With respect to this identification we have the following
matrix representations
\[
    \alp \ = \ \begin{pmatrix}
                    \id_L\\0
                \end{pmatrix}, \quad
    \bet \ = \ \begin{pmatrix}
                    0,& \id_N
            \end{pmatrix}, \quad
    \p^M \ = \ \begin{pmatrix}
        \p^L& f\\0&\p^N
    \end{pmatrix}
\]
where $f:N^*\to L^{*+1}$ is a morphism satisfying
\[
    f\circ\p^N \ = \ - \p^L\circ f.
\]
$f$ induces a morphism of the extended cohomology
$
    f_*: \calH^*(N) \ \longrightarrow \ \calH^{*+1}(L).
$
The long exact sequence of cohomology corresponding to \refe{exseq} has the
form
\begin{equation}\Label{E:longex}
 \begin{CD}
    \cdots @>{f_*}>> \H^i(L) @>{\alp_*}>> \H^i(M) @>{\bet_*}>>
    \H^i(N) @>{f_*}>> \H^{i+1}(L) @>{\alp_*}>>\cdots
 \end{CD}
\end{equation}\label{E:three}
Using this exact sequence one defines a natural isomorphism
\begin{equation}\Label{E:psi}
    \del: \det\H(L)\otimes \ \det\H(N) \ \overset{\sim}\longrightarrow \
    \det\H(M).
\end{equation}
Let $k_i$ and $c_i$ denote the kernel and cokernel of the morphism $f_\ast:
\H^{i-1}(N)\to \H^{i}(L)$, correspondingly. \refl{kerXcoker} gives a natural
isomorphism
\begin{eqnarray}\label{E:one}
\det \H^i(L)\otimes (\det\H^{i-1}(N))^\ast \, \stackrel\sim \to\,  \det
c_i\otimes (\det k_i)^\ast.
\end{eqnarray}
The isomorphism of \refp{seqofmod} applied to the exact
sequence\footnote{Formally \refp{seqofmod} deals with objects of the initial
category $\c$. However a similar statement with $M,M', M"\in \ob(\ec)$ is true
as follows from \refl{kerXcoker}} $$ 0\to c_i\to \H^i(M)\to k_{i+1}\to 0 $$
gives a natural isomorphism
\begin{eqnarray}\label{E:two}
\det c_i\otimes\det k_{i+1}  \stackrel\sim \to \det \H^i(M).
\end{eqnarray}
The alternating product of isomorphisms (\ref{E:one}) and (\ref{E:two}) produce
the isomorphisms
$$\det\H(L)\otimes \det\H(N) \to \det c\otimes (\det k)^\ast \to \det \H(M)$$
and their composition is denoted by $\delta$, see (\ref{E:psi}).

%
The following statement is the main result of this section.

\th{exseq=} Under the above conditions one has
\begin{equation}\Label{E:exseq=}
        \del\circ (\nu_L\otimes\nu_N) \circ \psi \ = \ \nu_M :\, \det M\to \det \H(M),
\end{equation}
where $\nu_L, \ \nu_M$, and $\nu_N$ are the canonical isomorphisms of the
determinant lines defined in \refp{detofcoh} and $\psi$ is the isomorphism of
\refp{seqofmod}. \eth

\begin{proof} This statement is similar to Theorem 3.2 of Milnor
\cite{Milnor66} stating that the torsion of an extension equals the product of
the individual torsions times the torsion of the exact homological sequence.
The full proof of Theorem \ref{E:exseq=} is quite long and technical; therefore
we have decided to publish it elsewhere. In the sequel this theorem is used
only once (in the proof of \reft{reltor}) and only in a special case when at
least one of the complexes $L, N, M$ is acyclic. In this special case
\reft{exseq=} admits a simpler proof which is still quite involved. We postpone
it to a further publication.
\end{proof}

\subsection{The relationship between the torsions}\Label{SS:reloftor}
{}We keep the notations introduced above. Fix $\c$-admissible scalar products
on $L$ and $N$. They determine a $\c$-admissible scalar product on $M$.
Although this scalar product on $M$ is not unique, its class in $\det M$depends
only on the classes in $\det L$ and $\det N$ represented by the given scalar
products on $L$ and on $N$. Let $\rho_L\in \det\H^\ast(L), \ \rho_M\in
\det\H^\ast(M)$, and $\rho_N\in \det\H^\ast(N)$  be the torsions of the
complexes $(L,\p^L)$, $(M,\p^M)$ and $(N,\p^N)$ respectively, cf.
\refd{torsion}. Let $\sig_L\in \det{L}, \ \sig_M\in \det{M}$, and $\sig_N\in
\det{N}$ be the elements induced by the scalar products, cf. \refss{projcomp}.
By definition, $
    \psi_{\alp,\bet}(\sig_M) \ = \ \sig_L\otimes \sig_N.
$
Hence we obain:
\th{torexseq} $\delta(\rho_L\otimes \rho_N) = \rho_M.$ \eth

\subsection{Torsion of the cone complex}\Label{SS:cone}
The above results can be reformulated in terms of the cone complex. Suppose
$(C,\p)$ and $(\tilC,\tilp)$ are chain complexes in $\c$. Let $f:C\to \tilC$ be
a chain morphism. Consider the cone of $f$, i.e., the chain complex $\Cone(f)$:
\begin{equation}\Label{E:cone}
 \begin{CD}
    0 \ \to\ C^0\oplus 0 @>{D_1}>> C^1\oplus\tilC^0 @>{D_2}>>
     \cdots @>{D_n}>> C^{n-1}\oplus \tilC^n@>{D_{n+1}}>> 0\oplus\tilC^n\ \to \ 0,
 \end{CD}
\end{equation}
where
\begin{equation}\Label{E:cone'}
    D_i \ = \   \begin{pmatrix}
                    -\p_{i}& 0\\ f_{i-1}& \tilp_{i-1}
                \end{pmatrix}.
\end{equation}

{}Fix $\c$-admissible scalar products on $C$ and $\tilC$. They induce a
$\c$-admissible scalar product on $\Cone(f)$ in an obvious way. Let $\rho_C\in
\det\calH(C)$, $\rho_{\tilC}\in \det\calH(\tilC)$ and $\rho_f\in
\det\H(\Cone(f))$ be the torsions of the complexes $(C,\p)$, $(\tilC,\tilp)$
and $\Cone(f)$ respectively, defined using these scalar products, cf.
\refd{torsion}.

{}From \reft{torexseq} we obtain:

\cor{torofcone}
Using the notation introduced above we have
\begin{equation}\Label{E:torofcone}
    \rho_f \ = \ \del_f\left(\rho_C\otimes(\rho_{\tilC})^\ast\right),
\end{equation}
where $(\rho_{\tilC})^\ast\in (\det(\H(\tilC)))^\ast$ denotes the dual of
$\rho_{\tilC}$. \ecor

This statement is equivalent to \reft{torexseq}. The dual $(\rho_\tilC)^\ast$
of the torsion appears because of the grading shift.

\section{Determinant class condition for chain complexes}

In this section we recall the determinant class condition for chain complexes
in Hilbertian categories. Historically, this condition appeared as a necessary
requirement for the $L^2$-torsion to be well defined. We show that the
determinant class condition for a chain complex $C$ is equivalent to the
$\tau$-triviality of the torsion parts of the extended $L^2$-cohomology
$\H(C)$. As we have shown in \S \ref{triviality}, the $\tau$-triviality implies
that the determinant lines of the torsion parts of the extended cohomology
could be canonically trivialized. Hence in this case the determinant line $\det
\H(C)$ can be canonically identified with the determinant line $\det H(C)$
where $H(C)$ is the reduced $L^2$-cohomology of $C$. This explains why under
the determinant class condition for a chain complex $C$ the torsion is well
defined as an element of $\det(H(C))$, as was proven in \cite{CarFarMat97}.

\defe{Dclasscomp} Let $\calC$ be a Hilbertian von Neumann category supplied
with a non-negative, normal, and faithful trace $\tau$. A chain complex
\begin{equation}\Label{E:Dclasscomp}
 \begin{CD}
    (C,\p): \ \
    0 \ \to \ C^0@>{\p_1}>> C^1@>{\p_2}>>\cdots @>{\p_{n-1}}>> C^n \ \to \ 0
 \end{CD}
\end{equation}
in $\c$ is said to be of $\tau$-determinant class if the induced maps
\begin{equation}\Label{E:Dclasscomp2}
    \p_j:\, C^{j-1}/\ker(\p_j) \ \to \ \cl(\im(\p_j))
\end{equation}
are $\tau$-isomorphisms for all $j=0\nek n$. \edefe

Compare \cite{BFKM} and \cite{CarFarMat97}.

Recall that the extended cohomology $\H^\ast(C)$ of the complex
\refe{Dclasscomp} is defined by formula (\ref{extended1}) and the reduced
cohomology is defined by formula \refe{reduced}.

\begin{lemma}{tautriv}
A chain complex $C$ as above is of $\tau$-determinant class if and only iff the
torsion parts $\calT(\H^i(C))$ of the extended cohomology are $\tau$-trivial
for all $i$.
\end{lemma}
\begin{proof}
The torsion part $\calT(\H^i(C))$ equals $(\p_j:\, C^{j-1}/\ker(\p_j) \to
\cl(\im(\p_j)))$, see \cite{Farber96}. Hence the result follows by comparing
Definition \ref{E:Dclasscomp} and \refd{tautriv}.

\end{proof}

\cor{ext=red} Given a chain complex \refe{Dclasscomp} of determinant class,
there is a canonical orientation preserving isomorphism between the determinant
lines of the extended and the reduced cohomology of $(C,\p)$:
\begin{equation}\Label{E:ext=red}\notag
    \det\calH^j(C) \ \overset{\sim}{\longrightarrow} \
    \det{H^j(C,\p)}, \quad j=0\nek n.
\end{equation}
\ecor
\prf Recall from \cite{Farber96} that the reduced cohomology is naturally
isomorphic to the projective part of the extended cohomology. Hence, in using
\refe{tor+proj} and \refe{directsum2}, we see that
$$\det \H^i(C) \simeq \det H^i(C)\otimes \det \calT(\H^i(C))\simeq \det
H^i(C)$$ assuming that $C$ is of $\tau$-determinant class. On the last stage we
have applied \refl{tautriv} and \refl{tautriv1} using the fact that the
determinant line $\det \calT(\H^i(C))$ of the torsion part is canonically
trivialized. \eprf

{}From Propositions~\ref{P:detofcoh} and \ref{C:ext=red} we obtain the
following:

\cor{detofredcoh} Given a chain complex \refe{Dclasscomp} of determinant class,
there is a canonical orientation preserving isomorphism
\begin{equation}\Label{E:detofredcoh}\notag
   \tilphi_C:\, \det{C} \ \overset{\sim}{\longrightarrow} \ \det{H^*(C,\p)}.
\end{equation}
\ecor

\subsection{The torsion of a determinant class complex}\Label{SS:tordetclass}
Using \refc{ext=red} we can view the torsion $\rho_C$ of a determinant class
complex $(C,\p)$ as a positive element of the determinant line
$\det(H^*(C,\p))$ of the reduced cohomology. In this case our notion of torsion
essentially reduces to the one considered in \cite{CarFarMat97}. More
precisely, if the complex $(C,\p)$ is of determinant class, then the element
\[
    \sig^*\otimes \psi_C(\rho_C) \ \in \ (\det{C})^*\otimes \det(H^*(C,\p))
\]
coincides with the torsion defined in \cite{CarFarMat97}.

\section{Combinatorial $L^2$-torsion  without the determinant class condition}\Label{S:combtor}

In this section we suggest a generalization of the classical construction of
the combinatorial torsion of Reidemeister, Franz and de Rham to the case of
infinite dimensional representations. Our torsion represents a nonzero element
of the determinant line of the {\it extended} $L^2$ cohomology. Intuitively,
the torsion is a density supported on the extended $L^2$ cohomology of the
polyhedron.

In more detail, given a finite polyhedron $K$ and a unimodular representation
of its fundamental group on a Hilbertian module $M$, the torsion invariant
defined in this section is a positive element of the determinant line $
    \det\calH^\ast(K;M).$
Here $\calH^\ast(K,M)$ denotes the extended $L^2$-homology of $K$ with
coefficients in $M$, as defined originally in \cite{Farber98}, \S 6.5. The
definition of $\H^\ast(K;M)$ is briefly repeated below.

The main advantage of the construction suggested here is that it requires no
additional assumptions. In particular, unlike the work of the previous authors
studying the construction of $L^2$ torsion, we do not require our complexes be
of determinant class.

We show that under the determinant class assumption the torsion defined in this
section coincides with the torsion defined in our previous work
\cite{CarFarMat97} where the torsion was understood as an element of the
determinant line of the {\it reduced} $L^2$ cohomology.

\subsection{Hilbertian bimodule}\Label{SS:bimodule}
Let $K$ be a finite cell complex. Denote by $\pi=\pi_1(K)$ its fundamental
group and by $C_\ast(\widetilde K)$ the cellular chain complex of the universal
covering $\widetilde K$ of $K$. The group $\pi$ acts on $C_\ast(\widetilde K)$
from the left and $C_\ast(\widetilde K)$ becomes a chain complex of free
finitely generated $\ZZ[\pi]$-modules with lifts of the cells of $K$
representing a free basis of $C_\ast(\widetilde K)$ over $\ZZ[\pi]$.

Let $M$ be a Hilbert space on which a von Neumann algebra $\A$ of bounded
linear operators acts from the right. As in \S \ref{SSS:hilbmod} we will denote
by $\B(M)=\B_{\A}(M)$ the commutant of $M$; recall that $\B_\A(M)$ is the ring
of all bounded linear maps $M\to M$ commuting with the right $\A$-action.

We shall assume that the algebra $\A$ is supplied with a finite, normal and
faithful trace $\tau$.

Consider a linear representation of the group $\pi$ in $M$. It is a group
homomorphism $\pi\to \B(M)\ =\ \B_{\A}(M)$. We think of $\pi$ as acting
linearly on $M$ from the left (via the representation $\pi\to \B(M)$); thus $M$
obtains an $(\pi-\A)$-bimodule structure. In this situation we say that $M$ is
an {\em Hilbertian $(\pi-\A)$-bimodule}.

A Hilbertian $(\pi-\A)$-bimodule $M$ is called {\em unimodular} if for every
element $g\in \pi$ the Fuglede-Kadison determinant $\Det_{\tau}(g) $ equals 1
where $g$ is viewed as an invertible linear operator $M\to M$ given by the
right multiplication by $g$.


\subsection{Flat Hilbertian bundle}\Label{SS:flbundle}
Any Hilbertian $(\pi-\A)$-bimodule $M$ determines a {\em flat Hilbertian bundle
over $K$ with fiber $M$}. Consider the Borel construction $
    \calE\ =\ {\widetilde K}\times_{\pi} M
$ together with the obvious projection map $\calE\to K$; it has a canonical
structure of a flat $\A$-Hilbertian bundle. Equivalently, it can be viewed as a
locally free sheaf of Hilbertian $\A$-modules.

\subsection{Examples}\Label{SS:exbimodule}
As the first example, consider the well-known situation when $\A={\calN}(\pi)$
is the von Neumann algebra of $\pi$ and $M=\ell^2(\pi)$ is the completion of
the group algebra of $\pi$ with respect to the canonical trace on it. Here $\A$
acts on $M=\ell^2(\pi)$ from the right and the group $\pi$ acts on $M$ from the
left. This $(\pi-\A)$-bimodule $M=\ell^2(\pi)$ is unitary and hence it is
unimodular, as well.

As a more general example consider the following. Let $V$ be a finite
dimensional unimodular representation of $\pi$. Let
$M=V\otimes_{\CC}\ell^2(\pi)$. Here the right action of $\A={\calN}(\pi)$ on
$M$ is the same as the action on $\ell^2(\pi)$ and the left action of $\pi$ is
the diagonal action: $
    g(x\otimes v)  =  gx\otimes gv
$ for $x\in \ell^2(\pi)$, $v\in V$, and $g\in \pi$.

\subsection{Construction of the torsion}\Label{SS:combtor}
Let $M$ be a $(\pi-\A)$-bimodule. Let
\begin{equation}\Label{E:cohain}
    C^\ast(K;M)\ =\ \Hom_{\ZZ[\pi]}(C_\ast(\widetilde K),M).
\end{equation}
It is a cochain complex of Hilbertian spaces (having right $\A$-action) and the
differentials are bounded linear maps commuting with the $\A$-action. More
precisely we may view complex (\ref{E:cohain}) as a cochain complex in the von
Neumann category $\calC=\H_\A$ of Hilbertian representations of $\A$, see \S
\ref{SSS:hilbmod}.

The cohomology $\H^\ast(K;M)$ of (\ref{E:cohain}) is an object of the extended
abelian category $\ec$. It is called {\it the extended cohomology} of $K$ with
coefficients in $M$. It was introduced originally in \cite{Farber95},
\cite{Farber96}.

By \refp{detofcoh} we obtain a natural isomorphism
\begin{equation}\Label{E:detC=detH}
 \begin{CD}
    \phi_{C^*(K;M)}: \, \det{C^\ast(K;M)}@>{\sim}>>\det{\H^\ast(K;M)}.
 \end{CD}
\end{equation}

Assume that the $(\pi-\A)$-bimodule $M$ is  unimodular. Then there is natural
isomorphism
\begin{equation}\Label{E:MchiK=C}
 \begin{CD}
    \psi_K:\, (\det{M})^{\chi(K)}@>{\sim}>> \det{C^\ast(K;M)},
 \end{CD}
\end{equation}
defined as follows. For any cell $e\subset K$ fix a lifting $\tilde e$ of $e$ in the
universal covering. Then the cells $\tilde e$ form a free $\ZZ[\pi]$-basis of the complex
$C_\ast(\widetilde K)$ and therefore induce an isomorphism
\begin{equation}\Label{E:M=C}
    \begin{CD}
        C^\ast(K;M) @>{\sim}>> \oplus_{e\in K} M
    \end{CD}
\end{equation}
between $C^\ast(K;M)$ and the direct sum of copies of $M$, one for each cell of
$K$. Thus, the determinant line $\det{C^\ast(K;M)}$ can be identified with
$(\det{M})^{\chi(K)}$ since the cells of odd dimension contribute negative
factors of $\det{M}$ into the total determinant line.

We only have to show that the above identification does not depend on the
choice of the liftings $\tilde e$. For any cell $e\subset K$ choose a group
element $g_e\in \pi$. Then the cells $g_e\tilde e$ form another set of
liftings. Consider the map $\oplus M \to \oplus M$ which is given by the
diagonal matrix with $g_e$ on the diagonal. The Fuglede-Kadison determinant of
this map equals 1 (since the representation $M$ is unimodular). Hence we see
that the isomorphism \refe{MchiK=C} is canonical.

\defe{combtor}
Fix a positive element $\sig\in \det{M}$. It defines an element
$\osig=\sig^{\chi(K)}\in (\det{M})^{-\chi(K)}.$ The combinatorial $L^2$ torsion
is defined by the formula
\begin{equation}\Label{E:combtor}
    \rho_K\ := \ \phi_{C^*(K,M)}\big(\psi_K(\osig)\big) \ \in\  \det\H^\ast(K;M).
\end{equation}
\edefe

In other words, the combinatorial $L^2$ torsion can be defined as follows:\,
let $\<\,,\,\>$ be an $\calC=\H_\A$-admissible scalar product on $M$ which
represents $\sig$, cf. \refss{equiv}. {}For each cell $e\in K$ fix a lifting
$\tile\in \tilK$. Then $\<\,,\,\>$ induces an admissible scalar product on
$C^\ast(K;M)$ via the isomorphism \refe{M=C}. The combinatorial torsion is the
torsion of the cochain complex $C^*(K;M)$ associated to this scalar product,
cf. \refd{torsion}.

\rem{combtor} 1. \ In the case when $\chi(K)=0$ the torsion $\rho_K$ is clearly
independent of the choice of $\sig$. This is always the case if $K$ is a closed
manifold of odd dimension.

2. \ If $\A=\CC$ we arrive at the classical definitions, cf.
\cite{Milnor61,Milnor62,Milnor66,Muller93,BisZh92}.

In the case when the cochain complex $C^\ast(K,M)$ is of determinant class,
\refd{combtor} reduces to the case considered in \cite{CarFarMat97}, cf.
\refss{tordetclass}. If, in addition, the reduced $L^2$-homology $H^\ast(K,M)$ vanishes,
we can identify the determinant line $\det{H^\ast(K,M)}$ with $\RR$ and so $\rho_K$ is
just a number. Under this assumption it was studied in \cite{CarMat92,Luck92,LuckRot91}.

3. Although our notation $\rho_K$ for the torsion invariant does not involve explicitly
the trace $\tau:\A\to \CC$, the whole construction (including the Fuglede-Kadison
determinants and the determinant lines) certainly depend on the choice of the trace
$\tau$, cf. Remark~4.5.3 of \cite{CarFarMat97}.
\erem

Recall that the classical Reidemeister-Franz torsion is not, in general, a homotopy
invariant, so one cannot expect homotopy invariance from our torsion invariant. But {\em
combinatorial invariance} holds in the following sense
\begin{Thm}[\textbf{Combinatorial Invariance}]\Label{T:combinv}
Let $K$ be a finite polyhedron and let $K^\prime$ be its subdivision. Suppose
that $M$ is a Hilbertian unimodular $(\pi-\A)$-bimodule. Let
\begin{equation}\Label{E:subdivH}
    \psi:\,\H^\ast(K^\prime;M) \ \to \ \H^\ast(K;M)
\end{equation}
be the isomorphism induced on the extended $L^2$ cohomology by the subdivision
chain map. Then the push-forward map
\begin{equation}\Label{E:subdivdet}
    \psi_\ast: \, \det\H^\ast(K';M)\ \to\ \det\H^\ast(K;M)
\end{equation}
preserves the torsions, i.e. it maps $\rho_{K'}$ onto $\rho_{K}$.
\end{Thm}
\prf It is enough to consider the elementary subdivision when a single
$q$-dimensio\-nal cell $e$ is divided into two $q$-dimensional cells $e_+$ and
$e_-$ introducing an additional separating $(q-1)$-dimensional cell $e_0$;
compare \cite{CarFarMat97}, proof of Theorem 4.6.

We have the exact sequence of free left $\ZZ[\pi]$-chain complexes
\begin{equation}\Label{E:seqofcomp}
 \begin{CD}
    0\to C_\ast(\widetilde K)@>{\psi}>>C_\ast(\widetilde K^\prime)\to D_\ast\to 0
 \end{CD}
\end{equation}
where the chain complex $D_\ast$ has nontrivial chains only in dimension $q$
and $q-1$ and $D_q$ and $D_{q-1}$ are both free of rank one. The free generator
of the module $D_q$ can be labelled with $e_+$ and the generator of $D_{q-1}$
can be labelled with the cell $e_0$ and then the boundary homomorphism is given
by $\partial(e_+)=e_0$. The exact sequence \refe{seqofcomp} induces the exact
sequence
\begin{equation}\Label{E:seqofcomp2}
    \begin{CD}
        0\to \Hom_{\ZZ[\pi]}(D_\ast;M)@>>> C^*(K';M) @>{\psi^*}>> C^*(K;M)\to 0.
    \end{CD}
\end{equation}
{}Fix an admissible scalar product $\<\,,\,\>$ on $M$. Using the cell
structures of $K$ and $K^\prime$, we then obtain canonically the scalar
products on the complexes $\Hom_{\ZZ[\pi]}(D_\ast,M)$, $C^*(K;M)$, $C^*(K';M)$.
The chain complex $\Hom_{\ZZ[\pi]}(D_\ast,M^*)$ is acyclic and so the canonical
isomorphism \refe{detofcoh} identifies the determinant line of $D_\ast$ with
$\RR$. The torsion of this complex equals $1$. \reft{combinv} now follows from
\refc{torofcone}. \eprf

\section{Extended cohomology of an elliptic complex}\Label{S:elcompl}

Throughout this section $\A$ denotes a finite von Neumann algebra with a fixed finite,
normal, and faithful trace $\tau$.

\subsection{Hilbertian $\A$-Bundles}\Label{SS:Abun}
Let $X$ be a closed manifold. A smooth bundle $p:\E\to X$ of topological vector
spaces (cf. Chapter~3 of \cite{Lang-DifMan}) is called a bundle of finitely
generated Hilbertian $\A$-modules or, simply, a Hilbertian $\A$-bundle if
\begin{enumerate}
\item $\E$ is equipped with a smooth fiberwise action $\rho: \E\times\A\to \E$
so that with this action each fiber $p^{-1}(x)$ is a finitely generated
Hilbertian $\A$-module, cf. \refsss{hilbmod}; \item there exist a finitely
generated Hilbertian $\A$-module $M$ such that $p:\E\to X$ is locally
isomorphic to the trivial bundle $p_0: X\times M\to X$ so that the local
isomorphism intertwines $p,p_0$ and the $\A$-action.
\end{enumerate}

Given a Hilbertian $\A$-bundle $\E$ over $X$ we denote by $C^\infty(X,\E)$ the space of
smooth sections of $\E$.

\subsection{Example: Flat Hilbertian $\A$-Bundles}\Label{SS:flatAbun}
Let $M$ be a finitely generated Hilbertian $(\pi-\calA)$-bimodule, cf.
\refss{bimodule}. Let $X$ be a connected, closed, smooth manifold with
fundamental group $\pi$ and let $\widetilde X$ denote the universal covering of
$X$. A flat Hilbertian $\A$-bundle with fiber $M$ over $X$ is an associated
bundle $p:{\E}\to X$, where $\E\ =\ (M\times\widetilde X)/\sim$ with its
natural projection onto $X$. Here $(m,x) \sim (gm,gx)$ for all $g\in \pi$,
$x\in \widetilde X$ and $m\in M$. Then $p:{\E}\to X$ is a a bundle of finitely
generated Hilbertian $\A$-modules.

\subsection{Elliptic complex of Hilbertian $\A$-modules}\Label{SS:elcomplex}
Let \/ $\E_0\nek \E_n$ \/ be Hilbertian $\A$-modules over a closed manifold $X$. Consider
a complex
\begin{equation}\Label{E:elcomplex}
 \begin{CD}
    0 \ \to \  C^\infty(X,\E_0) @>d_1>> C^\infty(X,\E_1) @>d_1>> \cdots @>d_{n}>>
    C^\infty(X,\E_n) \ \to \ 0,
 \end{CD}
\end{equation}
where
\begin{equation}\Label{E:Di}
    d_j:\, C^\infty(X,\E_{j-1}) \ \to \ C^\infty(X,\E_{j}),
    \qquad  j=1\nek n,
\end{equation}
are first order differential operators such that
 \(
    d_{j+1}\circ d_j  =  0.
 \)

The complex \refe{elcomplex} induces in a standard way the {\em symbol complex} over the
cotangent bundle $T^*X$ of $X$. The complex \refe{elcomplex} is called {\em elliptic} if
the symbol complex is exact outside of the zero section.

The complex \refe{elcomplex} is not a complex of Hilbertian spaces. Therefore, we can not
directly  define the the extended cohomology of this complex. However, following
\cite{Farber98}, we introduce in \refss{Sobolev} the Sobolev completion of
\refe{elcomplex}. We show that if the complex \refe{elcomplex} is elliptic then the
obtained extended cohomology does not depend on the choice of the Sobolev parameter $s$.

\subsection{Example: de Rham complex}\Label{SS:DeRcomp}
Given a flat Hilbertian $\A$-bundle $\E\to X$ over a closed connected manifold
$X$ one may consider the space of smooth differential $j$-forms on $X$ with
values in $\E$; this space will be denoted by $\Omega^j(X,\E)$. It is naturally
defined as a right $\A$-module. An element of $\Omega^j(X,\E)$ can be uniquely
represented as a $\pi$-invariant differential form on $\tilde X$ with values in
$M$, i.e. as a $\pi$-invariant element of $M\otimes_\CC \Omega^j(\widetilde
X)$. Here one considers the total (diagonal) $\pi$ action, that is, the tensor
product of the actions of $\pi$ on $\Omega^j(\widetilde X)$ and on $M$. More
precisely, if $\omega\in\Omega^j(\widetilde X)$ and $m\in M$, then
$m\otimes\omega$ is said to be $\pi$-invariant if $gm\otimes g^*\omega =
m\otimes\omega$ for all $g\in\pi$.

 A flat $\A$-linear connection on a flat Hilbertian $\A$-bundle $\E$ is
defined as an $\A$-homomor\-phism
\[
    \nabla= \n_{j+1} :\, \Omega^j(X,\E) \ \longrightarrow \ \Omega^{j+1}(X,\E)
\]
such that
\[
    \nabla(f\omega )\ = \ df\wedge\omega + f\nabla (\omega)
    \qquad\text{and}\qquad\nabla^2 = 0
\]
for any $\A$ valued smooth function $f$ on $X$ and for any $\omega\in
\Omega^j(X,\E)$. On a flat Hilbertian $\A$-bundle $\E$, as defined in
\refss{flatAbun}, there is a canonical flat $\A$-linear connection $\nabla$
which is given as follows: under the identification of $\Omega^j(X,\E)$ given
in the previous paragraph, one defines the connection $\nabla$ to be the de
Rham exterior derivative. If $\n$ is a flat connection on $\E$ then the de Rham
complex
\[
  \begin{CD}
    0 \ \longrightarrow \
    \Ome^0(X,\E) @>{\n_1}>> \Ome^1(X,\E) @>{\n_2}>>\cdots@>{\n_n}>> \Ome^n(X,\E) \
    \longrightarrow \ 0
  \end{CD}
\]
is an elliptic complex.

\subsection{The Laplacian}\Label{SS:lapl}
A Hermitian metric $h$ on a flat Hilbertian $\A$-bundle $p:\E\to X$ is a smooth
family of admissible scalar products on the fibers. Let us fix a smooth measure
$\mu$ on $X$. Then any Hermitian metric on $p:\E\to X$ defines the $L^2$ scalar
product
\begin{equation}\Label{E:L2scpr}
    \<s_1,s_2\> \ := \ \int_X\, h(s_1(x),s_2(x))\, d\mu(x), \qquad s_1,s_2\in C^\infty(X,\E),
\end{equation}
on $C^\infty(X,\E)$. We denote by $L^2(X,\E)$ the completion of the space
$C^\infty(X,\E)$ with respect to the scalar product (\ref{E:L2scpr}).

Suppose we are given an elliptic complex  \refe{elcomplex} of Hilbertian
$\A$-modules. Assume further that each bundle $\E_j, \, j=0\nek n$ is endowed
with a Hermitian metric $h_j$.

The {\em Laplacian} $\Delta_j$ is defined to be
\begin{equation}\Label{E:laplace}\notag
     \Delta_j \ = \ d_{j-1} d_{j-1}^{*} + d_{j}^* d_{j}:\,
     L^2(X,\E_j) \ \longrightarrow \ L^2(X,\E_j),
\end{equation}
where $d_j^{*}$ denotes the formal adjoint of \/ $d_j$ with respect to the $L^{2}$ scalar
product on \/ $L^2(X,\E_j)$. Note that, by definition, the Laplacian is a formally
self-adjoint operator which is densely defined. Since it is elliptic it is also
essentially self-adjoint. We denote by $\Delta_j$ the self adjoint extension of the
Laplacian.

\subsection{Spectral cut-off of an elliptic complex}\Label{SS:cutoff}
{}For every $\calI\subset \RR$ we denote by
\[
    L^2_\calI(X,\E_j) \ \subset  \ L^2(X,\E_j)
\]
the image of the spectral projection $P^j_\calI:L^2(X,\E_j)\to L^2(X,\E_j)$ of $\Del_j$
corresponding to $\calI$. The
following theorem was  proven by M. Shubin \cite[Th.~5.1]{Shubin98}%
\footnote{In \cite{Shubin98} the theorem is stated for the de Rham complex of a
flat Hilbertian $\A$-bundle. But the proof given there is valid for a general
elliptic complex.}.

\th{cutoff}
{}Fix $\eps>0$. Then
\begin{enumerate}
\item
$L^2_{[0,\eps]}(X,\E_j)\subset C^\infty(X,\E_j)$, i.e., $L^2_{[0,\eps]}(X,\E_j)$ consists
of smooth forms.
\item
The Hilbertian $\A$-module  $L^2_{[0,\eps]}(X,\E_j)$ is finitely generated, i.e., belongs
to the category $\calH^f_\calA$, cf. \refsss{hilbmod}.
\end{enumerate}
\eth

\subsection{The Sobolev complex}\Label{SS:Sobolev}
Let $W^s(X,\E_j)$ denote  the completion of the space of smooth sections in the topology
defined by the inner product
\begin{equation}\Label{E:sobprod}
    \<\omega_1,\omega_2\>_s = \big\<\,(I+\Delta_j)^{s}\omega_1,\omega_2\, \big\>,
    \qquad \ome_1,\ome_2\in C^\infty(X,\E_j).
\end{equation}

{}Fix $s>n$ and consider the following Sobolev complex
\begin{equation}\Label{E:Sobcomp}
 \begin{CD}
    0 \ \to \ W^s(X,\E_0)@>d_1>> W^{s-1}(X,\E_1)@>d_2>>\cdots
    \\ \cdots@>d_n>>W^{s-n}(X,\E_n) \ \to \ 0.
 \end{CD}
\end{equation}
This is a complex of Hilbertian $\A$-modules. Denote by $\calH^*_s(\E_*)$ the
extended cohomology of this complex.

Let
\[
    \Del_{j,s}:\,  W^{s-j}(X,\E_j) \ \to \  W^{s-j}(X,\E_j), \qquad j=0\nek n,
\]
denote the Laplacian of the complex \refe{Sobcomp}, defined as in \refe{laplace} but
using the scalar products \refe{sobprod}. Clearly,
\begin{equation}\Label{E:Del-Dels}
    \Del_{j,s} \ = \ (I+\Del_j)^{-1}\Del_j.
\end{equation}

Let $W^{s-j}_{[0,\del]}(X,\E_j)$ denote the image of the spectral projection of
$\Del_{j,s}$ corresponding to the interval $[0,\del]\subset \RR$. From \refe{Del-Dels} we
conclude
\begin{equation}\Label{E:A-As}
    L^2_{[0,\eps]}(X,\E_j) \ = \ W^{s-j}_{[0,\del]}(X,\E_j),
    \qquad\text{where}\quad \eps = \del(1+\del).
\end{equation}

\subsection{Extended cohomology of an elliptic complex}\Label{SS:extelliptic}
{}Fix $\del>0$ and set $\eps= \del(1+\del)$. By \reft{cutoff}, the subcomplex
\begin{equation}\Label{E:cutoffcomp}
 \begin{CD}
    0 \ \to \ L^2_{[0,\eps]}(X,\E_0)@>d_1>> L^2_{[0,\eps]}(X,\E_1)@>d_2>>\cdots
    @>d_n>>L^2_{[0,\eps]}(X,\E_n) \ \to \ 0
 \end{CD}
\end{equation}
of \refe{elcomplex} is a complex of finitely generated Hilbertian modules.

By \refe{A-As}, the inclusion
\begin{equation}\Label{E:inclusion}
    L^2_{[0,\eps]}(X,\E_j) = W^{s-j}_{[0,\del)]}(X,\E_j)
    \ \hookrightarrow \ W^{s-j}(X,\E_j)
\end{equation}
induces an isomorphism of the extended cohomology of the complexes
\refe{cutoffcomp} and \refe{Sobcomp}. Thus: {\it the cohomology of complexes
(\ref{E:Sobcomp}) and (\ref{E:cutoffcomp}) depends neither on the choice of
$\eps>0$ nor on the choice of the Sobolev parameter $s$}. We denote this
cohomology by $\H^*(\E_*)$ and refer to it as the {\em extended cohomology of
the elliptic complex \refe{elcomplex}}. Since $\H^*(\E_*)$ is isomorphic to the
extended cohomology of the finitely generated complex \refe{cutoffcomp}, it is
an object of the category $\E(\H^f_\A)$.

\subsection{Extended de Rham cohomology}\Label{SS:extDeRham}
Suppose that the elliptic complex \refe{elcomplex} is the de Rham complex of a
flat Hilbertian $\A$-bundle $\E$ with fiber $M$. Fix a Hermitian metric $h$ on
$\E$ and a Riemannian metric $g$ on $X$. These data define an $L^2$-scalar
product on $\Ome^j(X,\E)$ in a standard way.

We denote the extended cohomology of this complex by $\H^*(X,\E)$ and refer to
it as the {\em extended de Rham cohomology of $X$ with coefficients in $\E$}.

Let $W^{s}(\Ome^j(X,\E)), \ (j=0\nek n)$ denote the completion of the space
$\Ome^j(X,\E)$ of smooth forms with respect to the Sobolev scalar product
\refe{sobprod}. Then the complex \refe{Sobcomp} has the form
\begin{equation}\Label{E:SobDeRham}
 \begin{CD}
    0 \ \to \ W^s(\Ome^0(X,\E))@>\n>> W^{s-1}(\Ome^1(X,\E))@>\n>>\cdots
    \\ \cdots@>\n>>W^{s-n}(\Ome^n(X,\E)) \ \to \ 0.
 \end{CD}
\end{equation}

Let $K$ be a triangulation of $X$. Consider the cochain complex $C^*(K;M)$, cf.
\refe{cohain}. If $s>3n/2+1$, then the complex \refe{SobDeRham} consists of
continuous forms and the de Rham integration map
\[
    \tet:\,  W^{s-*}(\Ome^*(X,\E))  \ \longrightarrow \
    C^*(K;M)
\]
is defined. By the de Rham theorem for extended $L^2$-cohomology, cf.
\cite[\S7]{Farber98}, \cite{Shubin98}, this map is a homotopy equivalence. In
particular, it induces an isomorphism of the extended cohomology
\begin{equation}\Label{E:DeRhamTh}
    \tet_*:\, \calH^*(X;\E) \ \overset{\sim}\longrightarrow \ \calH^*(K;M).
\end{equation}

\section{Analytic $L^2$-torsion without the determinant class condition}\Label{S:DeRham}

We continue to use the notation introduced in the previous section. In
particular, $\A$ denotes a finite von Neumann algebra with a fixed finite,
normal, and faithful trace $\tau$.

\subsection{Torsion of the complex $L^2_{[0,\eps]}(X,\E_*)$}\Label{SS:toreps}
The $L^2$-scalar product \refe{L2scpr} on $C^\infty(X,\E_*)$ induces an admissible scalar
product on the complex $(L^2_{[0,\eps]}(X,\E_*),d_*)$ of finitely generated Hilbertian
$\A$-modules. We denote by
\begin{equation}\Label{E:smtor}
    \rho_{[0,\eps]} \ := \ \phi(\sig) \ \in \  \det\H^*(\E_*).
\end{equation}
its torsion as defined in \refd{torsion}.

\subsection{Torsion of the complex  $L^2_{(\eps,\infty)}(X,\E_*)$}\Label{SS:toreps2}
The complex $(L^2_{(\eps,\infty)}(X,\E_*),d_*)$ is acyclic and one would like
to define its torsion using the formula \refe{projcomp}. However, since there
is no trace on the category $\calH_\A$ we need to use the $\zet$-function
regularization of the determinant.

Let us recall the basic facts about the $\zet$-function of an elliptic operator on a
Hilbertian $\A$-bundle, cf. \cite{BFKM}. We will use the well known fact that for each
$\eps\ge 0$  the heat operators
\[
    e^{-t\Del_j|_{L^2_{(\eps,\infty)}(X,\E_j)}}, \qquad j=0\nek n,
\]
have smooth Schwartz kernels which are smooth sections of a bundle over
$X\times X$ with fiber the commutant of $M$, cf. \cite{BFKM,GromSh91,Luke72}.
The symbol $\Tr_\tau$ denotes application of the canonical trace on the
commutant (cf. \refss{trfgmod}) to the restriction of the kernels to the
diagonal followed by integration over the manifold $X$. The trace $\Tr_\tau$
vanishes on commutators of smoothing operators. See also
\cite{Mathai92,Lott92,GromSh91} for the case of the de Rham complex of a flat
bundle defined by the regular representation of the fundamental group.

It is shown in \cite[\S2.4]{BFKM} that the integral
\begin{equation}\Label{E:zet}
    \zet_{(\eps,\infty)}^j(s) \ = \
    \frac1{\Gam(s)}\, \int_{0}^\infty\, t^{s-1}\,
      \Tr_\tau(e^{-t\Del_j|_{L^2_{(\eps,\infty)}(X,\E_j)}})\, dt,
\end{equation}
converges for $\Re s<-\dim X/2$. Moreover, the function
$\zet_{(\eps,\infty)}^j(s)$ is analytic in $s$ for $\Re s<-\dim X/2$ and has a
meromorphic extension to the whole complex plane which is regular at $s=0$.

Set
\begin{equation}\Label{E:regdet}
    \log\Det \Del_j|_{L^2_{(\eps,\infty)}(X,\E_j)} \ = \
    \frac{d}{ds}|_{s=0}\, \zet_{(\eps,\infty)}^j(s).
\end{equation}

\lem{detxdet}
{}For every $\eps_1 >\eps_2>0$ we have
\begin{equation}\Label{E:detxdet}
     \Det\Del_j|_{L^2_{(\eps_2,\infty)}(X,\E_j)} \ = \
      \Det\Del_j|_{L^2_{(\eps_1,\infty)}(X,\E_j)}\cdot
        \Det_\tau \Del_j|_{L^2_{(\eps_2,\eps_1]}(X,\E_j)},
\end{equation}
where $\Det_\tau \Del_j|_{L^2_{(\eps_2,\eps_1]}(X,\E_j)}$ is the
Fuglede-Kadison determinant of the automorphism $\Del_j$ of the finitely
generated $\A$-module $L^2_{(\eps_2,\eps_1]}(X,\E_j)$. \elem \prf This easily
follows from the definitions using formula (\ref{E:projcomp}). \eprf

We now define the torsion $\rho_{(\eps,\infty)}$ using \refe{projcomp} and \refe{regdet},
i.e., we set
\begin{equation}\Label{E:tor>eps}
    \log\rho_{(\eps,\infty)} \ := \
    \frac12\, \sum_{j=0}^n\,
     (-1)^jj\frac{d}{ds}|_{s=0}\, \zet_{(\eps,\infty)}^j(s).
\end{equation}

{}From \refl{detxdet} we obtain:

\lem{epsinv} The product
\[
    \rho_{[0,\eps]}\cdot \rho_{(\eps,\infty)}
    \ \in \  \det\calH^*(\E_*)
\]
is independent of $\eps>0$.
\elem

\refl{epsinv} allows us to give the following definition.
\defe{analtor}
The analytic $L^2$ torsion of the complex $(C^\infty(X,\E_*),d_*)$ is defined
by the formula
\begin{equation}\Label{E:analtor}
    \rho_\E \ := \ \rho_{[0,\eps]}\cdot \rho_{(\eps,\infty)}
    \ \in \  \det\calH^*(\E_*).
\end{equation}
The  $L^2$-torsion of  the de Rham complex of a flat Hilbertian $\A$-bundle
$\E$, is called the analytic $L^2$-torsion of bundle $\E$. \edefe

Note that, in general, the torsion $\rho_\E$ depends on the choices of the
hermitian metrics on $\E_i$ and the measure on $X$. We shall show below in
\refs{chmu} that the de Rham torsion is independent of these choices if
$\dim{X}$ is odd and the bundle $\E$ is unimodular.

\section{The Cheeger-M\"uller type theorem}\Label{S:chmu}

In this section we recall the definition of the relative torsion of a flat
unimodular Hilbertian $\A$-bundle over a compact Riemannian manifold, see
\cite{BFK01}, \cite{CarMatMis99}. We show that the relative torsion equals to
an appropriately defined ratio of the de Rham and combinatorial torsions.
Burghelea, Friedlander and Kappeler \cite{BFK01} proved that the relative
torsion equals one for odd dimensional manifolds. Using this result we prove
here that for odd dimensional manifolds the analytic $L^2$ torsion coincides
with the combinatorial $L^2$ torsion, where both these torsions are viewed as
elements of the determinant line of the extended $L^2$ cohomology. Previously
this result was  proven in \cite{BFKM} under an additional assumption that the
Hilbertian bundle is of determinant class.

\subsection{The setting}\Label{SS:setreltor}
Let $X$ be a closed smooth manifold of dimension $\dim{X}=n$ with fundamental
group $\pi= \pi_1(X)$. Let $\calA$ be a finite von Neumann algebra and let $M$
be a finitely generated unimodular Hilbertian $(\pi-\A)$-bimodule, cf.
\refss{bimodule}.  Let $\calE$ be the flat Hilbertian $\calA$-bundle with fiber
$M$ over $X$, cf. \refss{flatAbun}. Fix a Hermitian metric $h$ on $\E$ (cf.
\refss{lapl}) and a Riemannian metric $g$ on $X$. Let
$
    \rho_\calE \  \in\ \det\calH^\ast(X;\E)
$
be the  analytic $L^2$-torsion, cf. \refd{analtor}.

{}Fix a cell decomposition $K$ of $X$ and a positive element of $\det{}M$. Let
$
    \rho_{K,M}\ \in\ \det\H^\ast(K;M)
$ denote the combinatorial $L^2$-torsion, cf. \refd{combtor}.

\subsection{The relative torsion}\Label{SS:reltor}
{}Fix a Sobolev index $s>n/2+1$. Then the Sobolev space $W^s(\Omega^*(X, \E))$
(cf. \refss{extDeRham}) consists of continuous forms and the de Rham
integration map
\[
    \tet:\, W^s(\Omega^*(X, \E)) \ \longrightarrow \ C^*(K; M)
\]
is defined. Consider the following sequence of chain maps:
\begin{equation}\Label{E:gsK}\notag
    (\Omega^*(X,\E),\nabla) \buildrel{(\Delta+I)^{-s/2}}\over {\longrightarrow}
    (W^s(\Omega^*(X,\E),\nabla) \buildrel{\theta_s}\over{\longrightarrow}
    (C^*(K; M),\p).
\end{equation}
The composition of these maps is denoted by $g(s,K)$. We form the cone complex:
\begin{equation}\Label{E:relcomp}\notag
    {\Cone}(g(s,K)) \ = \ C^*(K, M)\oplus\Omega^{*-1}(X, \E)
\end{equation}
with differential
\begin{equation}\Label{E:Drelcomp}
    D_i = D_i(s,K) \ = \ \begin{pmatrix}-\nabla_{i}&0\\g(s,K)_{i-1}&\p_{i-1}\end{pmatrix},
\end{equation}
where $\p$ is the coboundary map of the complex $C^*(K; M)$. It follows (as in
\cite[\S4]{CarMatMis99}, see also \cite[\S2]{BFK01}) that the cone complex is
acyclic. Hence the torsion of the cone complex is simply a positive real
number, which is defined as follows.

Set
\[
    \Del_j   \ :=  \ D_{j-1}D_{j-1}^*+D_{j}^*D_{j}.
\]
As in \refss{toreps2}, we can define the trace $\Tr_\tau$ of the heat kernel
$e^{-t\Del_j}$ and the $\zet$-function
\[
    \zet_{\Del_j}(\sig) \ = \ \frac1{\Gam(\sig)}\, \int_{0}^\infty\, t^{\sig-1}\,
      \Tr_\tau(e^{-t\Del_j})\, dt.
\]
It is shown in  \cite[\S6]{CarMatMis99} (cf. also \cite[\S2]{BFK01}) that this
function is defined and analytic for $\Re\sig>n/2$ and has a meromorphic
continuation to the whole complex plane which is regular at $\sig=0$.

\defe{reltor}
The relative torsion $\calR= \rho_{\Cone(g(s,K))}$ is defined by the formula
\begin{equation}\Label{E:reltor}\notag
    \log{\calR} \ = \ \frac{1}{2}\sum_j(-1)^jj\zeta'_{\Del_j}(0).
\end{equation}
\edefe

\subsection{Relations between the de Rham, the combinatorial,
   and the relative torsions}\Label{SS:reltor=}
The de Rham integration map \refe{DeRhamTh} induces an isomorphism
\begin{equation}\Label{E:DeRhamThDet}
    \tet_*:\, \det\H^\ast(X,E)
    \  \overset{\sim}{\longrightarrow} \ \det\H^\ast(K;M).
\end{equation}
We denote by
\begin{equation}\Label{E:DeRhamThDet2}
    \otet:\, \det\H^\ast(X,E)\,\otimes\, (\det\H^\ast(K;M))^*
    \  \overset{\sim}{\longrightarrow} \ \RR,
\end{equation}
the isomorphism induced by $\tet_*$.

Let
\[
    \rho_{K,M}^*\ \in\ (\det\H^\ast(K;M))^*
\]
be the dual of the combinatorial torsion $\rho_{K,M}$. It is the unique element
of $(\det\H^\ast(K;M))^*$ such that $\<\rho_{K,M}^*,\rho_{K,M}\>=1$; here
$\<\,,\,\>$ denotes the pairing between $(\det\H^\ast(K;M))^*$ and
$\det\H^\ast(K;M)$.

\th{reltor} One has
\begin{equation}\Label{E:reltor=}
     {\calR} \ = \ \otet\big(\rho_{\E}\otimes\rho_{K,M}^*\big).
\end{equation}
\eth

\prf Let $L^2\Ome^j(X,\E)$ denote the space of square integrable differential
forms on $X$ with values in $\E$. {}For every $\calI\subset \RR$ we denote by
$\Ome^j_{\calI}(X,\E)$ the intersection of the image of the spectral projection
$P^j_\calI:L^2\Ome^i(X,\E)\to L^2\Ome^j(X,\E)$ of the Laplacian
$\Del_j=\n_{i-1}\n_{i-1}^*+\n_i^*\n_i$ corresponding to $\calI$ with the space
of smooth forms. By \reft{cutoff} the image of the spectral projection
corresponding to $[0,\eps]$ consists of smooth forms. It follows that
$\Ome^j(X,\E)= \Ome^j_{[0,\eps]}(X,\E)\oplus \Ome^j_{(\eps,\infty)}(X,\E)$ for
every $\eps>0$.

Set
\[
    g_\calI(s,K)_j \ := \ g(s,K)_j|_{\Ome^j_{\calI}(X,\E)}.
\]
Then, for every $\eps>0$, we obtain an exact sequence of acyclic complexes
\begin{equation}\Label{E:decofrel}\notag
  \begin{CD}
    0 \ \to \ \Cone(g_{[0,\eps]}(s,K)) @>{j}>> \Cone(g(s,K)) @>\pi>>
        \Ome^*_{(\eps,\infty)}(X,\E) \ \to \ 0,
  \end{CD}
\end{equation}
where $j$ is the natural inclusion and $\pi$ is the natural projection. The
Carey-Mathai-Mishchenko lemma, cf. \cite{CarMatMis99}, \cite[Lemma~1.14]{BFK01}, implies
that
\begin{equation}\Label{E:decofreltor}
    \calR=\rho_{\Cone(g(s,K))} \ = \ \rho_{\Ome^*_{(\eps,\infty)}(X,\E)}\cdot
    \rho_{\Cone(g_{[0,\eps]}(s,K))}.
\end{equation}
{}From \refc{torofcone} we conclude that
\begin{equation}\Label{E:smallreltor}
    \rho_{\Cone(g_{[0,\eps]}(s,K))} \ = \
    \otet(\rho_{\Ome^*_{[0,\eps]}(X,\E)}\otimes \rho_{K,M}^*).
\end{equation}
Combining \refe{analtor}, \refe{decofreltor} and \refe{smallreltor} we obtain
\refe{reltor}. \eprf

\rem{reltor}  In the case when the extended cohomology $\H(X,\E)\simeq \H(K,M)$
is trivial and, hence, the torsions $\rho_\E$ and $\rho_{K,M}$ are just real
numbers, the right hand side of \refe{reltor=} is the ration of these numbers
and we obtain
\[
     {\calR} \ = \ {\rho_{\E}}/{\rho_{K,M}}.
\]

From \reft{reltor} it follows that the relative torsion is independent of the
Sobolev parameter $s$ (this was originally proven in \cite{CarMatMis99}). Also,
from \refe{reltor=} and from the invariance of the combinatorial torsion under
subdivisions we conclude that the relative torsion is independent of the cell
decomposition $K$ of $X$; this result is called the {\em combinatorial
invariance} of the relative torsion and was first proven in \cite{BFK01}. \erem

\subsection{The odd dimensional case and the Cheeger-M\"uller type theorem}\Label{SS:chmu}
Let us consider the case when the dimension $n$ of $X$ is odd. The Euler
characteristic of $K$ vanishes $\chi(K)=0$ and therefore the combinatorial
torsion $\rho_{K,M}$ is independent of the choice of the volume form on $M$,
cf. \refr{combtor}.1. Assume that the Hermitian metric $h$ on $\E$ is
unimodular. The main result of \cite{BFK01} (cf. Th.~0.1 in \cite{BFK01})
states that, in this situation the relative torsion equals one, i.e. $\calR \ =
\ 1.$ Combined with \reft{reltor} this result implies the following extension
of the Cheeger-M\"uller theorem

\th{chmu}
Suppose that the dimension of $X$ is odd and that the Hermitian metric $h$ on $\E$ is
unimodular. Then the isomorphism \refe{DeRhamThDet} identifies the analytic and the
combinatorial torsions
\begin{equation}\Label{E:chmu}
    \tet_*( \rho_\calE) \ = \ \rho_{K,M}.
\end{equation}
In particular, the analytic $L^2$ torsion $\rho_{\calE}$ does not depend on the
choice of the Riemannian metric on $X$ and the unimodular Hermitian metric on
$\E$. \eth 

Theorem~0.1 of \cite{BFK01} calculates the relative torsion in the general case
without assuming that the dimension of $X$ is odd and that the Hermitian metric
on $\E$ is unimodular. In this case $\log\calR$ is represented as an integral
over $X$ of some explicitly calculated locally defined differential form (this
form was originally found by Bismut and Zhang, \cite{BisZh92}, who considered
the case of finite dimensional bundle $\E$). We refer the reader to
\cite{BFK01} for details.

The fact that the relative torsion (and, hence, the analytic torsion) is
independent of the metrics is much easier to prove than the equality $\calR=1$.
Moreover, this fact is used in the proof of this equality, cf. \cite{BFK01}.
The independence of the metrics is proven in section~4 of \cite{BFK01}. Again,
the result of \cite{BFK01} is more general and calculates explicitly the
dependence of $\calR$ on the metrics in the case when the dimension of $X$ is
even (so called, {\em anomaly formula}). Using this result and \reft{reltor}
one can easily obtain a similar formula for the dependence of the analytic
torsion on the metrics. We leave the details for an interested reader.


\bibliographystyle{amsplain}
\providecommand{\bysame}{\leavevmode\hbox to3em{\hrulefill}\thinspace}
\providecommand{\MR}{\relax\ifhmode\unskip\space\fi MR }
\providecommand{\MRhref}[2]{%
  \href{http://www.ams.org/mathscinet-getitem?mr=#1}{#2}
} \providecommand{\href}[2]{#2}

\end{document}